\documentclass[aos,MSNbibl,nameyear,dvips]{arximspdf}
\usepackage{euscript}

\doi{10.1214/12-AOS1027} 
\volume{40}
\issue{5}
\pubyear{2012}
\firstpage{2421}
\lastpage{2451}

\makeatletter

\renewcommand{\tilde}{\widetilde}
\renewcommand{\hat}{\widehat}

\newcommand{\rrvert}{\vert}
\newcommand{\llvert}{\vert}

\newtheorem{theorem}{Theorem}
\newtheorem{lemma}{Lemma}
\newtheorem{definition}{Definition}

\newproclaim{remark}{Remark}

\newcommand{\iint}{\int\!\!\int}

\newcommand{\dddot}[1]{\hspace*{1.5pt}\dot{\vphantom{\Psi}}\hspace*{-0.2pt}\ddot{\hspace*{-1.3pt}#1}}

\renewcommand{\cite}{\citet}

\renewcommand{\a}{\alpha}
\renewcommand{\b}{\beta}
\newcommand{\Gam}{\Gamma}
\newcommand{\Del}{\Delta}
\newcommand{\eps}{\varepsilon}
\newcommand{\Lam}{\Lambda}
\newcommand{\lamb}{\lambda}
\newcommand{\sig}{\sigma}

\newcommand{\barlam}{\overline\lambda}

\newcommand{\NN}{{\mathbb N}}
\newcommand{\QQ}{{\mathbb Q}}
\newcommand{\RR}{{\mathbb R}}

\newcommand{\Cal}{\EuScript}

\newcommand{\cc}{{\Cal C}}
\renewcommand{\ll}{{\Cal L}}
\newcommand{\yy}{{\Cal Y}}

\newcommand{\trace}{\operatorname{trace}}

\newcommand{\abar}{\bar a}
\newcommand{\ghat}{\widehat g}
\newcommand{\that}{\widehat t_n}
\newcommand{\bhat}{\widehat b}
\newcommand{\BBhat}{\widehat{\mathbb{B}}}

\newcommand{\psidot}{\dot\psi}
\newcommand{\psiddot}{\ddot\psi}
\newcommand{\psidddot}{\dddot\psi}

\newcommand{\diag}{\operatorname{diag}}
\newcommand{\XXbar}{\overline{\XX}}
\newcommand{\ZZbar}{\overline{\ZZ}}
\newcommand{\zij}{z_{i,j}}
\newcommand{\zik}{z_{i,k}}
\newcommand{\etaik}{\eta_{i,k}}
\newcommand{\etaij}{\eta_{i,j}}
\newcommand{\lamiN}{\lambda_{i,N}}

\newcommand{\tz}{\widetilde z}

\newcommand{\tzik}{\widetilde z_{i,k}}

\newcommand{\tA}{\widetilde A}
\newcommand{\tB}{\widetilde B}
\newcommand{\tD}{\widetilde D}
\newcommand{\tdel}{\widetilde\delta}
\newcommand{\tH}{\widetilde H}
\newcommand{\tS}{\widetilde S}
\newcommand{\tll}{\widetilde\ll}
\newcommand{\txx}{\widetilde{\Cal X}}
\newcommand{\tyy}{\widetilde{\yy}}
\newcommand{\teta}{\widetilde\eta}
\newcommand{\tth}{\widetilde\theta}
\newcommand{\tgam}{\widetilde\gamma}
\newcommand{\txi}{\widetilde\xi}
\newcommand{\tlamiN}{\widetilde\lambda_{i,N}}
\newcommand{\tlam}{\widetilde\lambda}
\newcommand{\tphi}{\widetilde\phi}

\newcommand{\lldot}{\dot\ll}
\newcommand{\tlldot}{\dot{\widetilde\ll}}
\newcommand{\Hdot}{\dot H}
\newcommand{\Hddot}{\ddot H}

\newcommand{\tQQnaBN}{\widetilde{\QQ}_{n,a,\mathbb{B},N}}
\newcommand{\SPAN}{\operatorname{span}}

\newcommand{\Np}{{N_+}}
\newcommand{\sign}{\operatorname{sign}}

\newcommand{\rootlogn}{\sqrt{\log n}}

\newcommand{\scov}{\cc}

\newcommand{\var}{\operatorname{var}}
\newcommand{\id}{{d}}

\newcommand{\BB}{\mathbb{B}}
\newcommand{\PP}{\mathbb{P}}
\newcommand{\XX}{\mathbb{X}}
\newcommand{\ZZ}{\mathbb{Z}}

\newcommand{\ff}{\mathcal F}
\renewcommand{\tt}{\mathcal T}

\newcommand{\tmu}{\widetilde{\mu}}
\newcommand{\tK}{\widetilde{K}}

\newcommand{\PPnmuK}{\PP_{n,\mu,K}}
\newcommand{\QQnaB}{\QQ_{n,a,\mathbb{B}}}
\newcommand{\QQnaBN}{\QQ_{n,a,\mathbb{B},N}}

\newcommand{\hellinger}{\mathbf h}

\makeatother

\begin{document}
\begin{frontmatter}

\title{Estimation in functional regression for general exponential families}
\runtitle{Functional regression}

\begin{aug}
\author[A]{\fnms{Winston Wei} \snm{Dou}\thanksref{t3}\ead[label=e1]{Winston.Wei.Dou@aya.yale.edu}},
\author[A]{\fnms{David} \snm{Pollard}\thanksref{t2}\ead[label=e2]{David.Pollard@yale.edu}}
\and
\author[A]{\fnms{Harrison H.} \snm{Zhou}\corref{}\thanksref{t3}\ead[label=e3]{Huibin.Zhou@yale.edu}\ead[label=u1,url]{http://www.stat.yale.edu/}}
\runauthor{W. W. Dou, D. Pollard and H. H. Zhou}
\affiliation{Yale University}
\address[A]{Statistics Department\\
Yale University\\
24 Hillhouse Ave\\
New Haven, Connecticut\\
USA\\
\printead{e1}\\
\hphantom{E-mail: }\printead*{e2}\\
\hphantom{E-mail: }\printead*{e3}\\
\printead{u1}} 
\end{aug}

\thankstext{t2}{Supported in part by NSF Grant MSPA-MCS-0528412.}

\thankstext{t3}{Supported in part by NSF Career Award DMS-06-45676
and NSF FRG Grant DMS-08-54975.}

\received{\smonth{8} \syear{2011}}
\revised{\smonth{6} \syear{2012}}

%
\begin{abstract}
This paper studies a class of exponential family models whose
canonical parameters are specified as linear functionals of an
unknown infinite-dimensional slope function. The optimal minimax
rates of convergence for slope function estimation are established.
The estimators that achieve the optimal rates are constructed by
constrained maximum likelihood estimation with parameters whose
dimension grows with sample size. A change-of-measure argument,
inspired by Le Cam's theory of asymptotic equivalence, is used to
eliminate the bias caused by the nonlinearity of exponential family
models.
\end{abstract}

%

\end{frontmatter}

\section{Introduction} \label{intro}
There has been extensive exploratory and theoretical study of
\textit{functional data analysis} (FDA) over the past two decades. Two
monographs by Ramsay and Silverman (\citeyear{RamsaySilverman2002,RamsaySilverman2005})
provide comprehensive discussions on methods and applications.

Among many problems involving functional data, slope estimation in
functional linear regression has received substantial attention in
literature, for example, by \cite{CardotFerratySarda2003},
\cite{LiHsing2007} and \cite{HallHorowitz2007}. In particular,
\cite{HallHorowitz2007} established minimax rates of convergence
and proposed rate-optimal estimators based on spectral truncation
(regression on functional principal components). They showed that
the optimal rates depend on the smoothness of the slope function and
the decay rate of the eigenvalues of the covariance kernel of the
functional independent variable.

In this paper, we study optimal rates of convergence for slope
estimation in functional generalized linear models, for which little
theory is known. We introduce several new technical devices to
overcome the problems caused by nonlinearity of the link function.
To analyze our estimator, we establish a sharp approximation for
maximum likelihood estimators for exponential families parametrized
by linear functions of $N$-dimensional parameters, for an $N$ that
grows with sample size; see Lemma~\ref{AnBn}. We develop a
change-of-measure argument---inspired by ideas from Le Cam's theory
of asymptotic equivalence of models---to eliminate the effect of
bias terms caused by the nonlinearity of the link function; see
Sections~\ref{known} and~\ref{unknown}.

We consider problems where the observed data consist of independent,
identically distributed pairs $(y_i,\XX_i)$ where each $\XX_i$ is a
Gaussian process indexed by a compact subinterval of the real line,
which with no loss of generality we take to be $[0,1]$. Assume, for
each $i$, that the random variable $y_i$ conditional on the
process~$\XX_i$, follows a distribution $Q_{\lambda_i}$, where
$\{Q_\lambda\dvtx  \lambda\in\RR\}$ is a one-parameter exponential family.
The density function of $Q_{\lambda}$ is specified in equation
(\ref{gexpmodel}). We take parameter $\lambda_i$ to be a linear
functional of $\XX_i$ of the form
%
\begin{eqnarray}\label{lamidef}
\lambda_i &=& a +\int_0^1
\XX_i(t)\BB(t)\,dt\hspace*{80pt}
\nonumber\\[-8pt]\\[-8pt]
&&\eqntext{\mbox{for an unknown constant $a$ and an unknown $\BB\in
L^2[0,1]$.}}
\end{eqnarray}
Thus, the conditional joint distribution of $(y_1,\ldots, y_n)$
given $(\XX_1,\ldots, \XX_n)$ is the product measure
$\QQ_{n,a,\BB,\XX_1, \ldots, \XX_n} = \bigotimes_{i\leq n}
Q_{\lambda_i}$.
We abbreviate $\QQ_{n,a,\BB, \XX_1,\ldots, \XX_n}$ to
$\QQ_{n,a,\BB}$. Write $P_{\mu, K}$ for the distribution of each
$\XX_i$, where $\mu$ is the mean and $K$ is the covariance of
$\XX_i$. The joint distribution of the sample processes is then
$\PP_{n,\mu, K} = P_{\mu, K}^n$. Therefore, our models $\PP_{n,f}:=
\PP_{n,\mu,K}\QQ_{n,a,\BB}$, where $f=(K,a,\mu, \BB)$, are the joint
distributions of the sample $(y_1, \XX_1),\ldots, (y_n,\XX_n)$. The
parameter set $\ff\equiv\ff(R, \alpha, \beta)$ depend on universal
constants $R, \alpha$ and $\beta$. See Definition~\ref{ffdef} (in
Section~\ref{model}) for the precise specification of the parameter
set. The universal constant $\alpha$ controls the decay rate of
eigenvalues of kernel $K$, and the universal constant $\beta$
characterizes the ``smoothness'' of the slope function $\BB$.

Denote the corresponding norm and inner product in the space
$L^2[0,1]$ by $\|\cdot\|$ and $\langle \cdot,\cdot\rangle$. We focus
on the
estimation of $\BB$ using integrated squared error loss,
\[
L(\BBhat_n,\BB) = \|\BBhat_n -\BB\|^2 = \int
_0^1 \bigl(\BBhat_n(t) - \BB(t)
\bigr)^2 \,dt.
\]

The two main results are as follows.
%
\begin{theorem}[(Minimax upper bound)]\label{thmeupper} Under the
assumptions stated in Section~\ref{model}, there exists an
estimating sequence of $\BBhat_n$'s for which: for each $\eps>0$
there exists a finite constant $C_\eps$ such that
\[
\sup_{f\in\ff}\PP_{n,f} \bigl\{ \|\BBhat_n - \BB
\|^2>C_\eps n^{(1-2\beta)/(\alpha+2\beta)} \bigr\}<\eps\qquad\mbox{for all
large enough $n$.}\vadjust{\goodbreak}
\]
\end{theorem}
%
\begin{theorem}[(Minimax lower bound)]\label{thmelower} Under the
assumptions stated in Section~\ref{model},
\[
\liminf_{n\to\infty}n^{(2\beta-1)/(\alpha+2\beta)}\sup_{f\in\ff
}\PP_{n,f} \|
\BBhat_n - \BB\|^2 >0 \qquad\mbox{for every estimator $\{
\BBhat_n\}$.}
\]
\end{theorem}

Two closely related works in the area of functional generalized
linear models are \cite{MullerStadtmuller2005} and
\cite{CardotSarda2005}, which provided theory for the convergence
rate in functional generalized linear models. However, the rate
optimality was not studied. In addition,
\cite{MullerStadtmuller2005} established an upper bound for rates
of convergence assuming the negligibility of the bias due to the
approximation of the infinite-dimensional model by a sequence of
finite-dimensional models, the issue we overcome by
using a change-of-measure argument. 
By contrast, more theoretical
advances have been achieved in the functional linear regression setting,
not only for estimation but also for
prediction. For example, \cite{CaiHall2006} and
\cite{CrambesKneipSarda2009} derived optimal rates of convergence
for prediction in the fixed and random design cases. See also,
\cite{CardotMasSarda2007} which derived a CLT for prediction in the
fixed and random design cases and \cite{CardotJohannes2010} which
established a minimax optimal result for prediction at a random
design using thresholding estimators. In a companion study to our
paper, Dou [(\citeyear{Dou2010}), Chapter 5] considers optimal prediction in
functional generalized linear regressions with an application to the
economic problem of predicting recessions from the U.S. Treasury
yield curve.

Our minimax upper bound result (Theorem~\ref{thmeupper}) is proved
in Section~\ref{upperproof}. The minimax lower bound result
(Theorem~\ref{thmelower}) is established in
Section~\ref{lowerproof}. The proof of Theorem~\ref{thmeupper}
depends on an approximation result (Lemma~\ref{AnBn}) for maximum
likelihood estimators in exponential family models for parameters
whose dimensions change with sample size. As an aid to the reader,
we present our proof of Theorem~\ref{thmeupper} in two stages. In
Section~\ref{known}, we assume that both the mean $\mu$ and the
covariance kernel $K$ are known. This allows us to emphasize the key
ideas in our proofs. We proceed in Section~\ref{unknown} to the case
where $\mu$ and $K$ are estimated.
The proofs for the lemmas are collected together in
Section~\ref{lemmaproof}. Some of them invoke the
perturbation-theoretic results collected in the supplemental
article [\cite{DouPollardZhou2012Supplemental}].

\section{Regularity conditions} \label{model}
Let $\{Q_\lambda\dvtx  \lambda\in\RR\}$ be a one-parameter exponential family,
%
\begin{equation}
\label{gexpmodel}\id Q_\lambda/ \id Q_0 =
f_{\lambda}(y):= \exp\bigl(\lambda y -\psi(\lambda)\bigr)\qquad\mbox{for all }
\lambda\in\RR.
\end{equation}
Necessarily\vspace*{1pt} $\psi(0)=0$.
Remember that $e^{\psi(\lambda)} = Q_0 e^{\lambda y}$ and that the
distribution $Q_\lambda$ has mean $\psidot(\lambda)$ and
variance $\psiddot(\lambda)$.
\begin{remark*}
We may assume that $\psiddot(\lambda)>0$ for every real $\lambda$.
Otherwise we would have $0=\psiddot(\lamb_0)=\var_{\lambda_0}(y)=Q_0
f_{\lambda_0}(y)(y-\psidot(\lambda_0))^2$ for some $\lambda_0$,
which would
make $y=\psidot(\lambda_0)$ for almost all $y$ under $Q_0$ and hence
$Q_\lambda\equiv
Q_{\lambda_0}$ for every~$\lambda$.\vspace*{-2pt}
\end{remark*}
\begin{remark*}
The main results in this paper can be extended to the functional
quasi-likelihood regression models [see, e.g., \cite{Wedderburn74}]
as follows:
\[
y_i = \mu_i + \sigma_i
\varepsilon_i,
\]
where
\[
\mu_i = g \biggl(a + \int_\tt\BB(t)
\XX_i(t)\,dt \biggr) \quad\mbox{and}\quad \sigma_i = v (
\mu_i )\qquad \mbox{with known } g \mbox{ and } v.
\]
However, a
main goal of this paper is to provide a better understanding of the
difficulties caused by nonlinearity in functional data analysis
models and to propose a general approach to tackle them. The
exponential families can provide a good representation of the
quasi-likelihood regression models to this end. One of the gains of
specifying exponential families is to simplify the proofs while
still achieving our main goal and covering the most broadly used
models, such as the functional logistic regression model, the
functional probit regression model and the functional poisson
regression model. The general nonparametric setting where the link
functions $g$ and $v$ are unknown is studied by
\cite{MullerStadtmuller2005}, assuming the negligibility of the
bias due to the approximation of the infinite-dimensional model by a
sequence of finite-dimensional models. Without ignoring the bias,
the problem becomes much more difficult and would be an interesting
topic for future research.\vspace*{-2pt}
\end{remark*}
\begin{remark*}
A natural extension of our model is the classical generalized linear
model with nuisance parameters $\phi$ as follows:
\[
y_i | \XX_i \sim f_{\lambda_i, \phi}(y) \qquad\mbox{with }
\lambda_i = a + \int_\tt\BB(t)
\XX_i(t)\,\id t
\]
and
\[
f_{\lambda, \phi}(y):= \exp\bigl[\alpha_1(\phi) \bigl(
\lambda y -\psi(\lambda)\bigr) + \alpha_2(\phi, y) \bigr],
\]
where $\alpha_1(\phi) > 0$ so that for each $\phi\in
\RR^d$ we have an exponential family. Under some regularity
conditions on the known functions $\alpha_1(\cdot)$ and
$\alpha_2(\cdot)$, the exact maximum likelihood estimation analysis
and the lower bound argument of this paper can still be employed to
derive minimax results for the slightly more general setting.\vspace*{-2pt}
\end{remark*}

We assume:
\begin{longlist}[{${(\dddot\Psi)}$}]
\item[{$({\ddot\Psi})$}]
For each $\eps>0$ there exists a finite constant $C_\eps$ for
which $\psiddot(\lambda)\le C_\eps\exp(\eps\lambda^2)$ for
all $\lambda\in\RR$. Equivalently, $\psiddot(\lambda)\le
\exp(o(\lambda^2) )$ as $|\lambda|\to\infty$.\vadjust{\goodbreak}

\item[{${(\dddot\Psi)}$}]
There exists an increasing real function $G$ on $\RR^+$ such that
\[
\bigl|\psidddot(\lambda+h)\bigr|\le\psiddot(\lambda)G\bigl(|h|\bigr)
\qquad\mbox{for all $\lambda$ and $h$}.
\]
Without loss of generality we assume $G(0)\ge1$.
\end{longlist}



We also assume the observed data are i.i.d. pairs $(y_i,\XX_i)$ for
$i=1,\ldots,n$, where:
\begin{longlist}[{({X})}]
\item[{({X})}]
Each $\{\XX_i(t)\dvtx 0\le t\le1\}$ is distributed like $\{\XX(t)\dvtx 0\le
t\le
1\}$, a Gaussian process with mean $\mu(t)$ and covariance kernel $K(s,t)$.


\item[{({Y})}]
$y_i|\XX_i\sim Q_{\lambda_i}$ with $\lambda_i=a+\langle\XX_i,\BB
\rangle$
for an unknown $\{\BB(t)\dvtx 0\le t\le1\}$ in $L^2[0,1]$ and $a\in\RR$.
\end{longlist}
%
\begin{definition}\label{ffdef} For real constants $\a>1$ and
$\b>(\a+3)/2$ and $R>1$, define $\ff= \ff(R,\a,\b)$ as the set of
all $f=(K, a,\mu,\BB)$ that satisfy the following conditions:
\begin{longlist}[{({K})}]
\item[{({K})}]
The covariance kernel is square integrable with respect to Lebesgue measure
and has an eigenfunction expansion (as a compact operator on $L^2[0,1]$)
\[
K(s,t) = \sum_{k\in\NN}\theta_k
\phi_k(s)\phi_k(t),
\]
where the eigenvalues $\theta_k$ are decreasing with
$
Rk^{-\a}\ge\theta_k \ge\theta_{k+1}+ (\a/R) k^{-\a-1}
$.

\item[{({a})}]\label{a} $|a|\le R$.

\item[{$({\mu})$}]
$\|\mu\|\le R$.

\item[{({B})}]
$\BB$ has an expansion
$
\BB(t) = \sum_{k\in\NN}b_k \phi_k(t)
$
with $|b_k|\le Rk^{-\b}$, for the eigenfunctions defined by the kernel $K$.
\end{longlist}
\end{definition}
\begin{remark*} The purpose of this paper is not to offer a universally optimal
estimation procedure, but to provide a theory for the principal
components regression in nonlinear models of functional data. As in
\cite{HallHorowitz2007} and \cite{CaiHall2006}, among others,
assumptions {({K})} and {({B})} set up a natural
theoretical framework to
justify and analyze the principal components regression. In
practice, principal components analysis has been one of the most
widely and successfully used statistical methods. One example of
successful application of principal components analysis is in
analyzing the relationship between U.S. Treasury zero-coupon yield
curves, which is a typical functional data, and the macroeconomic
activities [see, e.g., \cite
{Dou2010,EstrellaHardouvelis91,Wright2006}]. In this analysis,
the fixed basis such as wavelet basis or fourier basis fails to give
a sparse representation of the yield curve data. Admittedly, under
different regularity assumptions, by design the principal components
regression approach may not be applicable, and accordingly, other
estimation methods such as wavelet basis or spline basis may have
better performance; see, for example, \cite
{EfromovichKoltchinskii2001,CrambesKneipSarda2009}. In \cite
{EfromovichKoltchinskii2001}, the authors
discussed an approach of using two different bases, one is for the
slope function and the other is for the covariance kernel operator.
This technique can be applied to some cases where the principal
components regression fails. Nevertheless, the results in
Efromovich and Koltchinskii
[(\citeyear{EfromovichKoltchinskii2001}), Theorem~3.1] requires a lower
level of noise in the covariance kernel and a higher degree of smoothness
of the slope function in order to allow tractability in more
severely ill-posed settings.

\end{remark*}
\begin{remark*} The awkward lower bound for $\theta_k$ in
assumption {({K})} implies,
for all $k<j$,
%
\begin{equation}\label{incrth}
\theta_k -\theta_j\ge R^{-1}\int
_k^j \a x^{-\a
-1}\,dx = R^{-1}
\bigl(k^{-\a}-j^{-\a} \bigr).
\end{equation}
If $K$
and $\mu$ were known, we would only need the lower bound $\theta_k\ge
R^{-1}k^{-\a}$ and not the lower bound for $\theta_k-\theta_{k+1}$. As
explained by Hall and Horowitz [(\citeyear{HallHorowitz2007}), page 76], the stronger
assumption is needed when one estimates the individual
eigenfunctions of $K$. Note that the subset of $L^2[0,1]$ in
which $\BB$ lies, denoted as $\mathcal B_K$, depends on $K$. We regard the
need for the stronger assumption on the eigenvalues and the irksome
assumption {({B})} as artifacts of the method of proof, but
we have
not yet succeeded in removing either assumption.
\end{remark*}
\begin{remark*}
We discuss two extreme cases to help understand the regularity
assumption $\beta> (\alpha+ 3)/2$. One case is that the
eigenvalues $\{\theta_k\}$ decay exponentially fast and the slope
coefficients $\{b_k\}$ decay with polynomial rates, where
essentially we have $\alpha$ is much larger than $\beta$, for which
it can be shown that the optimal rate of convergence is just
logarithmic. The other case is that the eigenvalues $\{\theta_k\}$
decay polynomially fast, and the slope coefficients $\{b_k\}$ decay
with exponential rates, where essentially we have $\beta$ is much
larger than $\alpha$, for which it can be shown that the optimal
convergence rate is nearly parametric up to a logarithmic term.
\end{remark*}



\section{Methodology}\label{method}
In this section we introduce the methodology to construct a sequence
of estimators, which achieve the optimal rates of convergence stated
in Theorem~\ref{thmeupper}. Our estimation features a two-step
procedure. We first truncate at the first $N$ principal components
and replace the original model $\PP_{n,f}$ by the truncated
model $\widetilde{\PP}_{n,f,N}$ defined in (\ref{tPnfN}). The
choice of $N$ depends on an estimation-approximation trade-off:
oversized $N$ can compromise the performance of the MLE maximizing
(\ref{bhathat}), whereas undersized $N$ can make the model
misspecification between $\PP_{n,f}$ and its finite-dimensional
approximation $\widetilde{\PP}_{n,f,N}$\vadjust{\goodbreak} nonnegligible. Second, we
further truncate the MLE at $m < N$ to form our estimator in
(\ref{bhat}). The choice of $m$ depends on the standard
variance-bias tradeoff as in nonparametric estimation problems. See
Section~\ref{known} for more details.

Under the assumptions {({X})} and {({K})} from
Section~\ref{model}, the
process $\XX_i$ admits the eigen decomposition
\[
\XX_i (t) - \mu(t) =\ZZ_i(t) = \sum
_{k\in\NN}\zik\phi_k(t).
\]
The random variables $\zik:= \langle\ZZ_i,\phi_k\rangle$ are independent
with $\zik\sim N(0,\theta_k)$.

Because $\mu$ and $K$ are unknown, we estimate them in the usual
way:
%
\begin{equation}
\label{tmu} \tmu(t)=\XXbar(t) = n^{-1}\sum
_{i\le n} \XX_i(t)
\end{equation}
and
%
\begin{eqnarray}\label{tK}
\tK(s,t) &=& (n-1)^{-1}\sum_{i\le n} \bigl(
\XX_i(s)-\XXbar(s) \bigr) \bigl(\XX_i(t)-\XXbar(t)
\bigr)
\nonumber\\[-8pt]\\[-8pt]
&=& (n-1)^{-1}\sum_{i\le n} \bigl(
\ZZ_i(s)-\ZZbar(s) \bigr) \bigl(\ZZ_i(t)- \ZZbar(t)
\bigr),
\nonumber
\end{eqnarray}
which has spectral representation
%
\begin{equation}
\label{tdecomp} \tK(s,t) = \sum_{k\in\NN}
\tth_k \tphi_k(s)\tphi_k(t)
\end{equation}
with $\tth_1\geq
\tth_2\geq\cdots\geq\tth_{n-1}\geq0$. In fact we must
have $\tth_k=0$ for $k\ge n$ because all the
eigenfunctions $\tphi_k$ corresponding to nonzero $\tth_k$'s must
lie in the $(n-1)$-dimensional space spanned by $\{\ZZ_i-\ZZbar\dvtx
i=1,2,\ldots,n\}$.

Using the first $N$ [to be specified in (\ref{N})] principal
components, we can approximate the original infinite-dimensional
model $\PP_{n,f}$ by a sequence of truncated finite-dimensional
models
%
\begin{equation}
\label{tPnfN}\widetilde{\PP}_{n,f,N} = P_{n,\mu,
K}\widetilde{
\QQ}_{n,a,\BB, N, \XX_1, \ldots, \XX_n},
\end{equation}
where
$\widetilde{\QQ}_{n,a, \BB, N, \XX_1, \ldots, \XX_n}:=
\bigotimes_{i\leq n}Q_{\widetilde{\lambda}_{i,N}}$ with $y_i |
\XX_1,\ldots, \XX_n \sim Q_{\widetilde{\lambda}_{i,N}}$ and
%
\begin{equation}
\label{tlamN}\widetilde{\lambda}_{i,N} = \widetilde{b}_0 +
\sum_{1\leq k\leq N}\widetilde{b}_k (
\widetilde{z}_{i,k} - \widetilde{z}_{\cdot k} ),
\end{equation}
where $\widetilde{b}_0 = a +
\langle\BB,\XXbar\rangle$, and $\widetilde{b}_k =
\langle\BB,\widetilde{\phi}_k\rangle$ for $k\geq1$, and
$\widetilde{z}_{i,k} = \langle\ZZ_i,\widetilde{\phi}_k\rangle$
for all
$i,k$, and $\widetilde{z}_{\cdot k} = n^{-1}\sum_{i\leq
n}\widetilde{z}_{i,k} = \langle\ZZbar,\widetilde{\phi}_k\rangle$. And
hence $\widetilde{z}_{i,k} - \widetilde{z}_{\cdot k} = \langle\ZZ_i
- \ZZbar,\widetilde{\phi}_k\rangle= \langle\XX_i - \XXbar
,\widetilde{\phi}_k\rangle$. We
abbreviate $\widetilde{\QQ}_{n,a, \BB, N, \XX_1, \ldots, \XX_n}$
to $\widetilde{\QQ}_{n,a, \BB, N}$ in the rest of the paper. We
introduce the following matrices
and vectors for the purpose of notational convenience. Define:
\begin{itemize}
\item$z_i:= (z_{i,1},\ldots, z_{i, N})'$ and $\tz_i:= (\tz
_{i,1},\ldots,\tz_{i,N})'$;\vspace*{1pt}
\item$z_{\cdot}:= (z_{\cdot1}, \ldots, z_{\cdot N})'$ and $\tz_{\cdot
}:= (\tz_{\cdot1 },\ldots,\tz_{\cdot N })'$;\vadjust{\goodbreak}
\item$D:= \diag(1,\theta_1, \ldots, \theta_N)^{1/2}$, where
$\theta_k$'s are
defined in assumption {({K})};\vspace*{1pt}
\item$\tD:= \diag(1,\tth_1,\ldots,\tth_N)^{1/2}$, where $\tth_k$'s
are defined in (\ref{tdecomp});\vspace*{1pt}
\item$\xi_i:= (1, z_i')'$ and $\txi_i:= (1,\tz_i'-{\widetilde
z_{\cdot}}')'$;\vspace*{1pt}
\item$\eta_i:= D^{-1}\xi_i$ and $\teta_i:= D^{-1}\txi_i$;\vspace*{1pt}
\item$\gamma:= (b_0, b_1,\ldots, b_N)'$ and $\tgam:=(\widetilde
b_0,\widetilde b_1,\ldots,\widetilde b_N)'$.
\end{itemize}
Thus, equation (\ref{tlamN}) can be rewritten as
%
\begin{equation}
\label{tlamN2} \widetilde{\lambda}_{i,N} = \txi_i'
\tgam= \teta_i'D\tgam.
\end{equation}

We estimate $\BB$ by
%
\begin{equation}
\label{bhat}\widehat{\BB}_n (t) = \sum_{k\leq
m}
\widehat{b}_k\widetilde{\phi}_k (t),
\end{equation}
where $(\widehat{b}_0,
\ldots, \widehat{b}_N)$ is the conditional MLE for the truncated
model $\widetilde{\PP}_{n,f,N}$, and $m$ is the optimal cutoff point
according to the
variance-bias tradeoff with $m < N$. More precisely, $(\widehat{b}_0,
\ldots,
\widehat{b}_N)$ is chosen to maximize the following conditional (on
the $\XX_i$'s) log likelihood over $g\equiv(g_0,g_1,\ldots, g_N)'$
in $\RR^{N+1}$:
%
\begin{equation}
\label{bhathat} \ll_n(g) = \sum_{i\leq n}y_i
\bigl(\txi_i'g \bigr) - \psi\bigl(
\txi_i'g \bigr)
\end{equation}
%
with cutoff points $m$ and $N$ chosen as
%
\begin{equation}
\label{m}
m\asymp n^{1/(\alpha+2\beta)}
\end{equation}
and
%
\begin{equation}
\label{N} N\asymp n^\zeta\qquad\mbox{with }(2+2\a)^{-1} > \zeta> (
\a+2\b-1)^{-1}.
\end{equation}
Note that $N$ is much
larger than $m$. Such a $\zeta$ exists because the assumptions
$\a>1$ and $\b>(\a+3)/2$ imply $\a+2\b-1>2+2\a$. The universal
constants $\alpha$ and $\beta$
characterize the decay rate of the eigenvalues of kernel $K$ and the
smoothness of slope function $\BB$ defined
in Definition~\ref{ffdef}.

\section{\texorpdfstring{Proof of Theorem \protect\ref{thmeupper}}{Proof of Theorem 1}}\label{upperproof}
The proof of Theorem~\ref{thmeupper} is divided into two stages.
In the first stage, we prove the theorem assuming that the mean
$\mu$ and the covariance kernel $K$ are known. This case is
relatively simple and of course artificial, but it captures the
essence of our proof. For the Gaussian case, this is reduced to the
setting considered in \cite{GoldenshlugerTsybakov2001}. In the
second\vspace*{1pt} stage where $\mu$ and $K$ are unknown, we show that using the
natural estimates $\widetilde{\mu}$ and $\widetilde{K}$ as in
(\ref{tmu}) and (\ref{tK}) will not affect the result achieved in the
first stage. 

In Section~\ref{majorlemma} we state the technical lemmas which
serve as building blocks for establishing the main theorems. Their
proofs are postponed to the Section~\ref{lemmaproof}. In
Section~\ref{known} we prove Theorem~\ref{thmeupper} assuming
$\mu$ and $K$ are known, and then in Section~\ref{unknown} we apply
Lemma~\ref{txxn} to complete the proof of Theorem~\ref{thmeupper}
with unknown $\mu$ and $K$.

\subsection{Technical lemmas}\label{majorlemma}
We write the lemmas in a notation that makes the applications in
Sections~\ref{known} and~\ref{unknown} more straightforward. The
notational cost is that the parameters are indexed by
$\{0,1,\ldots,N\}$ in Lemmas~\ref{AnBn} and~\ref{AB}. Each of
the lemmas stated in this subsection is a general result. 

We first introduce an approximation result for maximum likelihood
estimators in exponential family models for parameters whose
dimensions change with sample size. This lemma combines ideas
from \cite{Portnoy88} and from \cite{HjortPollard93}.
For each square matrix $A$, its spectral norm is defined by its
largest absolute value of the eigenvalues, that is, $\|A\|_2:=
\sup_{|v|\leq1}|Av|$ where $|v|$ denotes the $l^2$ norm of vector
$v$. The proof can be found in Section~\ref{proofAnBn}.
%
\begin{lemma}\label{AnBn}
Let $Q_\lambda$ be the one-parameter exponential family distribution
defined as in (\ref{gexpmodel}) and satisfying regularity
condition {${(\dddot\Psi)}$}. Suppose $\xi_1,\ldots,\xi_n$
are (nonrandom)
vectors in $\RR^{N+1}$. Suppose $\QQ=\bigotimes_{i\le n}Q_{\lambda_i}$
with $\lambda_i=\xi_i'\gamma$ for a
fixed $\gamma=(\gamma_0,\gamma_1,\ldots,\gamma_{N})'$ in $\RR^{N+1}$.
Under $\QQ$, the coordinate maps $y_1,\ldots,y_n$ are independent
random variables with $y_i\sim Q_{\lambda_i}$.

The log-likelihood for fitting the model is
\[
L_n(g)=\sum_{i\le n}\bigl(
\xi_i'g\bigr) y_i -\psi\bigl(
\xi_i'g\bigr)\qquad\mbox{for $g\in\RR^{N+1}$},
\]
which is maximized (over $\RR^{N+1}$) at the MLE $\ghat$
$(\mbox{$=$}\ghat_n)$. Define $\eta_i:= D^{-1}\xi_i$ for some nonsingular
matrix $D$, and define the matrix
\[
J_n:= \sum_{i\leq n}\xi_i
\xi_i'\psiddot(\lambda_i) =
nDA_nD' \qquad\mbox{with } A_n:= \frac1n\sum
_{i\le
n}\eta_i\eta_i'
\psiddot(\lambda_i).
\]
Assume $B_n$ is another nonsingular matrix for which
%
\begin{equation}
\label{assumptionA} \|A_n-B_n\|_2\le\bigl(2
\bigl\|B_n^{-1}\bigr\|_2\bigr)^{-1}
\end{equation}
and assume
%
\begin{equation}\label{maxeta}
{\max_{i\le n}}|\eta_i| \le\frac{\eps\sqrt{n}/(N+1)}{ G(1) \sqrt{32\|
B_n^{-1}\|_2}} \qquad\mbox{for some
$0<\eps<1$},
\end{equation}
where $G(\cdot)$ is defined
as in regularity condition {${(\dddot\Psi)}$}. Then, for each
set of vectors
$\kappa= \{\kappa_0,\ldots,\kappa_M\}$ in $\RR^{N+1}$ there is a set
$\yy_{\kappa,\eps}$ with $\QQ\yy_{\kappa,\eps}^c<2\eps$ on which
\[
\sum_{0\le j\le M}\bigl|\kappa_j'(
\ghat-\gamma)\bigr|^2\le\frac{
6\|B_n^{-1}\|_2 } {n\eps} \sum
_{0\le j\le M}\bigl|D^{-1}\kappa_j\bigr|^2.
\]
\end{lemma}
\begin{remark*}
This is a quite general result. In this paper, we are interested in
one particular case where $\kappa_j$ have all elements equal to zero
except the $j$th element that equals one and $D =
\diag(\theta_0,\ldots, \theta_M)$. In this case, the result can be
rewritten as
\[
\sum_{0\leq j\leq M}(\widehat{g}_j -
\gamma_j)^2 \leq\frac{
6\|B_n^{-1}\|_2 } {n\eps}\sum
_{0\leq j\leq M}\theta_j^{-2}.
\]
\end{remark*}
The following approximation result for random matrices will be
invoked in order to apply Lemma~\ref{AnBn} to show Theorem
\ref{thmeupper}. The proof can be found in Section~\ref{proofAB}.
%
\begin{lemma}
\label{AB} 
Suppose
$\{\eta_{i,k}\dvtx i,k\geq1\}$ are i.i.d. standard normal random
variables. Let
%
\begin{equation}
\label{A}A_n:= n^{-1} \sum_{i\le
n}
\eta_i\eta_i'\psiddot\bigl(
\gamma'D\eta_i\bigr),
\end{equation}
where $\gamma= (\gamma_0,
\gamma_1, \ldots, \gamma_N)', \eta_i =
(1,\eta_{i,1},\ldots,\eta_{i,N})' $ and $D = \operatorname{diag}(D_0,
\ldots, D_N)$.
Define $B_n:= \PP A_n$, and assume $\psi$ satisfies regularity
condition {$({\ddot\Psi})$}. If we have $ \sum_{k\geq
1}D_k^2\gamma_k^2<\infty$ and $N=o (n^{1/2} )$, then it
follows that
\[
\bigl\|B_n^{-1}\bigr\|_2=O(1) \quad\mbox{and}\quad \PP
\|A_n-B_n\|_2^2=o(1).
\]
\end{lemma}

The following lemma establishes a bound on the Hellinger distance
between members of an exponential family, which plays a key role in our
change-of-measure argument. We write $\hellinger(\cdot,\cdot)$ for
the Hellinger distance. If both $P$ and $Q$ are dominated by some
measure $\nu$, with densities $p$ and $q$, then $
\hellinger^2(P,Q):=\nu(\sqrt{p} - \sqrt{q})^2$. The proof can be found
in Section~\ref{hellinger}.
%
\begin{lemma}\label{expfactslemma} 
Suppose $\{Q_\lambda\dvtx  \lambda\in\RR\}$ is an exponential family
defined as in (\ref{gexpmodel}) and satisfies regularity condition
{${(\dddot\Psi)}$}. Then
\[
\hellinger^2(Q_\lambda,Q_{\lambda+\delta})\le
\delta^2\psiddot(\lambda) \bigl(1 +|\delta| \bigr)G\bigl(|\delta|\bigr)\qquad \forall
\lambda, \delta\in\RR.
\]
Here $G(\cdot)$ is defined in the
condition {${(\dddot\Psi)}$}.
\end{lemma}

The following lemma provides a maximal inequality for
weighted-chi-square variables, which easily leads to maximal
inequalities for Gaussian processes and multivariate normal vectors.
These inequalities will be repeatedly invoked. The proof can be found
in Section
\ref{proofgaussian}.
%
\begin{lemma}\label{gaussian} 
Suppose $\{\eta_{i,k}\dvtx i,k\geq1\}$ are i.i.d. standard normal random
variables. Let
\[
W_i=\sum_{k\in\NN}\tau_{i,k}
\eta_{i,k}^2 \qquad\mbox{for $i=1,\ldots,n$.}
\]
If the $\tau_{i,k}$'s are nonnegative constants
with $T:=\max_{i\le n}\sum_{k\in\NN}\tau_{i,k} < \infty$, then
it follows that
\[
\PP\Bigl\{\max_{i\le n} W_i > 4T(\log n + x)\Bigr\} <
2e^{-x} \qquad\mbox{for each $x\ge0$.}
\]
\end{lemma}

The following lemma is to guarantee that the estimation of $\BB$
using $\widetilde{\mu}$ and $\widetilde{K}$ basically has the same
accuracy as using $\mu$ and $K$. We need some terminology before
formally introducing the lemma, and these notations introduced below
apply to the rest of the paper. When we want to indicate that a
bound involving constants $c$, $C$, $C_1,\ldots$ holds uniformly over
all models indexed by a set of parameters~$\ff$, we write $c(\ff)$,
$C(\ff)$, $C_1(\ff),\ldots\,$. By the usual convention for eliminating
subscripts, the values of the constants might change from one
paragraph to the next: a constant $C_1(\ff)$ in one place need not be
the same as a constant $C_1(\ff)$ in another place. For sequences of
constants $c_n$ that might depend on $f\in\ff$, we write
$c_n=O_\ff(1)$ and $o_\ff(1)$ and so on to show that the asymptotic
bounds hold uniformly over $\ff$. Denote $H_p$ and $\widetilde{H}_p$
to be orthogonal projection operators associated with
$\SPAN\{\phi_1,\ldots, \phi_p\}$ and
$\SPAN\{\widetilde{\phi}_1,\ldots, \widetilde{\phi}_p\}$,
respectively,\vspace*{1pt} where $\phi_k$'s are the eigenfunctions defined in
assumption (K), and $\widetilde{\phi}_k$'s are their sample
approximations defined in (\ref{tdecomp}). We also need to define the
following key quantities:
\begin{itemize}
\item
$\tS:=\diag(\sig_0,\ldots,\sig_N)$ with $\sig_0=1$ and
$\sig_k=\sign(\langle\phi_k,\tphi_k\rangle)$ for $k\ge1$.\vspace*{1pt}

\item
$\Del:=\tK-K$, where $\tK$ is defined in (\ref{tK}).\vspace*{1pt}









\item
$\tA_n:= n^{-1}\sum_{i\le
n}\teta_i\teta_i'\psiddot(\tlam_{i,N})$, where $\teta_i$ and
$\tlam_{i,N}$ are defined in Section~\ref{method}.\vspace*{1pt}

\item
$\tB_n:= \widetilde{S} B_n \widetilde{S}$, where $B_n$ is defined
in (\ref{B}).
\end{itemize}
The proof of Lemma~\ref{txxn} can be found in Section~\ref{prooftxxn}.
%
\begin{lemma}\label{txxn} 
Assume the regularity conditions in Section~\ref{model} hold. Let $m$
and $N$ are chosen according to
(\ref{m}) and (\ref{N}), respectively. For each $\eps>0$ there exists a
set $\txx_{\eps,n}$, depending
on $\mu$ and $K$, with
\[
\sup_\ff\PPnmuK\txx_{\eps,n}^c<\eps\qquad\mbox{for all
large enough $n$}
\]
and on which, for some constant $C_\eps$ that does not depend
on $\mu$ or $K$:
{\renewcommand\thelonglist{(\roman{longlist})}
\renewcommand\labellonglist{\thelonglist}
\begin{longlist}
\item\label{txxnDel}
$\|\Del\|\le C_\eps n^{-1/2}$;

\item\label{txxnZZ}
$\max_{i\le n}\|\ZZ_i\| \le C_\eps\rootlogn$ and
$\|\ZZbar\|\le C_\eps n^{-1/2}$;\vspace*{1pt}

\item\label{txxnHm}
$\|(\tH_m-H_m)\BB\|^2 = o_\ff(n^{(1-2\beta)/(\alpha+
2\beta)} )$;\vspace*{1pt}

\item\label{txxnHN}
$\|(\tH_N-H_N)\BB\|^2 = O_\ff(n^{-1-\nu} )$ for some
universal constant $\nu>0$;

\item\label{txxnteta}
$ \max_{i\le n}|\teta_i|^2 = o_\ff(\sqrt{n}/N)$;\vspace*{2pt}

\item\label{txxntA}
$\|\tS\tA_n \tS- A_n\|_2 = o_\ff(1)$.
%
\end{longlist}}
\end{lemma}
%

\subsection{\texorpdfstring{Proof of Theorem \protect\ref{thmeupper} with known Gaussian distribution}
{Proof of Theorem 1 with known Gaussian distribution}}\label{known}
Initially we suppose that $\mu$ and $K$
are known. We emphasize that this simpler case serves as an
intermediate step to the more interesting unknown distribution case,
and it captures the essential idea of the proof of
Theorem~\ref{thmeupper}. 

Remember under $\QQ_{n,a,\BB}$, the $y_i$'s are independent,
conditional on $\XX_1,\ldots,\break \XX_n$, with $y_i\sim Q_{\lambda_i}$ and
\[
\lambda_i = a +\langle\XX_i,\BB\rangle=
b_0 + \sum_{k\in\NN
}\zik b_k\qquad
\mbox{where $b_0=a+\langle \mu,\BB\rangle$.}
\]
Our task is to estimate the $b_k$'s with sufficient accuracy so that
we are able to estimate $\BB(t) = \sum_{k\in\NN}b_k\phi_k(t)$
within an error of order $n^{(1-2\b)/(\a+2\b)}$. In fact it
will suffice to estimate the component $H_m\BB$ of $\BB$ in the
subspace spanned by $\{\phi_1,\ldots,\phi_m\}$ with $m\asymp
n^{1/(\a+2\b)}$ because
%
\begin{equation}\label{tailBB}
\bigl\|H_m^\perp\BB\bigr\|^2 = \sum
_{k>m}b_k^2 = O_\ff
\bigl(m^{1-2\b}\bigr) =O_\ff\bigl(n^{(1-2\beta
)/(\alpha
+ 2\beta)} \bigr).
\end{equation}

One might try to estimate the coefficients $(b_0,\ldots,b_m)$ by
choosing $\ghat=(\ghat_0,\ldots,\ghat_m)'$ to maximize a conditional
log likelihood over all $g = (g_0, g_1,\ldots, \break g_m)'$
in $\RR^{m+1}$:
\[
\ll_{n,m}(g):=\sum_{i\le n} y_i
\biggl(g_0 + \sum_{1\le k\le
m}\zik
g_k \biggr)-\psi\biggl(g_0 + \sum
_{1\le k\le m}\zik g_k \biggr).
\]
To this end one might try to appeal to Lemma~\ref{AnBn} stated at
the beginning of the previous subsection, with $\kappa_j$ equal to
the unit vector with a $1$ in its $j$th position for $j\le m$ and
$\kappa_j=0$ otherwise. That would give a bound for $\sum_{k\le
m}(\ghat_k-b_k)^2$. Unfortunately, we cannot directly invoke the
lemma with $N=m$ to estimate $\gamma' =(b_0,b_1,\ldots,b_N)$ when we
replace $\QQ$, $D$, $\xi_i$ and $\eta_i$ (notations) in Lemma
\ref{AnBn} by $\QQ_{n,a,\BB}$ (defined in Section~\ref{intro}), $D$,
$\xi_i$ and $\eta_i$ (defined in Section~\ref{method}),
respectively,
because $\lambda_i\ne\xi_i'\gamma$, a bias problem.
\begin{remark*} We could modify Lemma~\ref{AnBn} to allow $\lambda_i=
\xi_i'\gamma+\mathrm{bias}_i$, for a suitably small bias term, but at
the cost of extra regularity conditions and a more delicate
argument. The same difficulty arises whenever one investigates the
asymptotics of maximum likelihood estimators with the true
distribution outside the model family, that is, MLE under model
misspecification.
\end{remark*}

Instead, we use a two-stage estimation procedure,
%
\begin{equation}
\label{truebhat}\widehat{\BB}_n = \sum_{k\leq
m}
\widehat{b}_k\phi_k,
\end{equation}
where\vspace*{1pt} $(\widehat{b}_0, \ldots,
\widehat{b}_N)$ is the conditional MLE for the truncated model and
$m\leq N$. More precisely, $(\widehat{b}_0, \ldots, \widehat{b}_N)$
is chosen to\vadjust{\goodbreak} maximize the following conditional (on the $\XX_i$'s)
log likelihood over $g\equiv(g_0,g_1,\ldots, g_N)$ in $\RR^{N+1}$:
\[
\ll_{n,N}(g):=\sum_{i\leq
n}y_i
\bigl(\xi_i'g \bigr) - \psi\bigl(\xi_i'g
\bigr)
\]
with cutoff
points $m$ and $N$ chosen as in (\ref{m}) and (\ref{N}),
respectively. Note that this estimator differs from that in
(\ref{bhat}) in the sense that it uses $\phi_k$ and $z_{i,k}$ instead
of the approximation correspondences $\widetilde{\phi}_k$ and
$\widetilde{z}_{i,k} - \widetilde{z}_{\cdot k}$. This two-stage
estimation procedure eliminates the bias term by a change-of-measure
argument conditional on the $\XX_i$'s. We present the proof in the
following three steps.


%
%

\subsubsection*{Step 1}
From the analysis above, one can see that the key in our proof is
the change-of-measure argument and the application of Lemma
\ref{AnBn}. In this step, we construct a high probability set such
that for each realization of the $\XX_i$'s on the set the
assumptions of Lemma~\ref{AnBn} are satisfied.

Define $\gamma$, $\xi_i$, $D$ and $\eta_i$ as in Section
\ref{method}. Note that in this case $\eta_{i,j} =
z_{i,j}/\sqrt{\theta_i}$ for all $i,j$, and hence the $\eta_{i,j}$'s
are i.i.d. standard normal variables. We define matrix $A_n$ as in
(\ref{A}),
%
\begin{equation}
\label{B} A_n = n^{-1} \sum_{i\le
n}
\eta_i\eta_i'\psiddot\bigl(
\gamma'D\eta_i\bigr) \quad\mbox{and}\quad B_n:=
\PPnmuK A_n.
\end{equation}
Now, let us define
${\Cal X}_n={\Cal X}_{\ZZ,n}\cap{\Cal X}_{\eta,n}\cap{\Cal
X}_{A,n}$, where
%
\begin{eqnarray}
\label{xxnZZ}
{\Cal X}_{\ZZ,n} &:=& \Bigl\{{\max_{i\le n}}\|\ZZ_i
\|^2 \le C_0\log n\Bigr\},
\\
\label{xxneta}
{\Cal X}_{\eta,n} &:=& \Bigl\{{\max_{i\le n}}|\eta_i|^2
\le C_0 N\log n\Bigr\},
\\
\label{xxnA}
{\Cal X}_{A,n} &:=& \bigl\{\|A_n-B_n
\|_2\le\bigl(2\bigl\|B_n^{-1}\bigr\|_2
\bigr)^{-1}\bigr\}.
\end{eqnarray}
If we choose a large enough universal constant $C_0=C_0(\ff)$,
Lemma~\ref{gaussian} ensures that $\PPnmuK{\Cal X}_{\ZZ,n}^c\le2/n$ and
${\PPnmuK{\Cal X}_{\eta,n}^c\le2/n}$ by choosing ${\tau_{i,k} =
\theta_i}$ and ${\tau_{i,k} = \{i\leq N\}}$, respectively, for all
$i,k$; and Lemma~\ref{AB} shows that
\[
\bigl\|B_n^{-1}\bigr\|_2=O_\ff(1)\quad\mbox{and}\quad
\PPnmuK\|A_n-B_n\|_2^2=o_\ff(1),
\]
thus $\PPnmuK{\Cal X}_{A,n}^c = o_\ff(1)$. And hence,
%
\begin{equation}
\label{xn}\PPnmuK{\Cal X}_{n}^c\leq\PPnmuK{\Cal
X}_{\ZZ,n}^c + \PPnmuK{\Cal X}_{\eta,n}^c
+ \PPnmuK{\Cal X}_{A,n}^c = o_\ff(1).
\end{equation}

\subsubsection*{Step 2} In the previous step, we show the
assumptions of Lemma~\ref{AnBn} are satisfied on the set ${\Cal X}_n$. In
this step, we show that the change-of-measure argument is ready to
work. Let us consider the truncated model
\[
\QQnaBN:= \bigotimes_{i\le n}Q_{\lamiN} \qquad\mbox{with $\lamiN:=
\xi_i'\gamma$.}
\]
``Change of measure'' means to view the data $y_1,\ldots, y_n$ as if
they are generated from the conditional joint distribution $\QQ_{n,
a, \BB, N}$, though the true distribution is $\QQ_{n,a,\BB}$. In
this step, we show that the divergence caused by replacing $\QQnaB$
by $\QQnaBN$ is small enough that it will not
compromise the asymptotic results. 
A~common control of this divergence is the total variation distance
between $\QQnaBN$ and $\QQnaB$. We show that there exists a sequence
of nonnegative constants $c_n$ of order $o_\ff(\log n)$ such that
%
\begin{equation}\label{TVbound}
\|\QQ_{n,a,\BB}-\QQnaBN\|_{\mathrm{TV}}^2 \le
e^{2c_n}\sum_{i\le
n}|\lambda_i-
\lamiN|^2 \qquad\mbox{on ${\Cal X}_n$.}
\end{equation}
To
establish inequality (\ref{TVbound}) we use the bound
\[
\|\QQ_{n,a,\BB}-\QQnaBN\|_{\mathrm{TV}}^2\le
\hellinger^2(\QQ_{n,a,\BB},\QQnaBN) \le\sum
_{i\le
n}\hellinger^2(Q_{\lambda_i},Q_{\lamiN}).
\]
By Lemma~\ref{expfactslemma}
\[
\hellinger^2(Q_{\lambda_i},Q_{\lamiN}) \le
\delta_i^2\psiddot(\lambda_i) \bigl(1+|
\delta_i| \bigr)G\bigl(|\delta_i|\bigr) \qquad\mbox{with }
\delta_i:= \lambda_i - \lamiN,
\]
where
%
\begin{eqnarray}
\label{deltai} |\delta_i| &=&|\lambda_i-\lamiN| = \bigl|
\langle\ZZ_i,\BB\rangle-\langle H_N\ZZ_i,
\BB\rangle\bigr| =\bigl|\bigl\langle\ZZ_i,H_N^\perp\BB
\bigr\rangle\bigr|
\nonumber\\[-8pt]\\[-8pt]
&\le& \|\ZZ_i\|\bigl\|H_N^\perp\BB\bigr\| \le
O_\ff\bigl(\sqrt{N^{1-2\b}\log n } \bigr) =
o_\ff(1).
\nonumber
\end{eqnarray}
Because $\delta_i = o_\ff(1)$ for each $i$, we know all the
$ (1+|\delta_i| )G(|\delta_i|)$ factors can be bounded by a
single $O_\ff(1)$ term.

Further, for $(a,\BB,\mu,K)\in\ff(R,\a,\b)$ and with the
$\|\ZZ_i\|$'s on the set ${\Cal X}_n$,
%
\begin{equation}
\label{lambdai} |\lambda_i| \le|a| +\bigl(\|\mu\|+\|\ZZ_i\|\bigr)
\|\BB\| \le C_2\sqrt{\log n}
\end{equation}
for some constant $C_2=C_2(\ff)$. Assumption {$({\ddot\Psi
})$} then ensures
that all the $\psiddot(\lambda_i)$ are bounded by a single
$\exp(o_\ff(\log n) )$ term.

Therefore, inequality (\ref{TVbound}) is proved to hold. This
bound for total variation distance legitimates the change-of-measure
argument in the next step.

\subsubsection*{Step 3} We apply the change-of-measure argument and
Lemma~\ref{AnBn} to complete the
proof.
On the set ${\Cal X}_n$, we can apply Lemma~\ref{AnBn} directly with
$\QQ=
\QQnaBN$, because the conditions of Lemma~\ref{AnBn} hold: inequality
(\ref{assumptionA}) holds by construction of ${\Cal X}_n$ and inequality
(\ref{maxeta}) holds for large enough $n$ because
\[
{\max_{i\le n}}|\eta_i|^2 \le O_\ff(N
\log n) = o_\ff(\sqrt{n}/N).
\]
In the equation above, the first inequality is due to the construction
of ${\Cal X}_n$, and the second equality is due to $N =
o_\ff(n^{1/(2+2\alpha)} )$.

For each realization of the $\XX_i$'s lying in ${\Cal X}_n$, we invoke
Lemma~\ref{AnBn}, with $\eta_i$, $A_n$, $B_n$, $D$ and $\QQ$
(notations) in Lemma~\ref{AnBn} replaced by $\eta_i$, $A_n$, $B_n$,
$D$ and $\QQnaBN$ defined in this subsection, respectively, and it
gives a high probability set $\yy_{m,\eps}$ with
$\QQnaBN\yy_{m,\eps}^c<2\eps$ on which
\[
\sum_{1\le k\le m}|\bhat_k-b_k|^2
= O_\ff\biggl(n^{-1}\sum_{1\le
k\le m}
\theta_k^{-1} \biggr) = O_\ff
\bigl(m^{1+\a}/n\bigr)=O_\ff\bigl(n^{(1-2\beta
)/(\alpha+ 2\beta)} \bigr),
\]
which implies
\[
\|\BBhat_n-\BB\|^2 = \sum_{1\le k\le m}|
\bhat_k-b_k|^2 + \sum
_{k>m}b_k^2 = O_\ff
\bigl(n^{(1-2\beta)/(\alpha+ 2\beta)} \bigr).
\]

From inequality (\ref{TVbound}) it follows, for a large enough
constant $C_\eps$, that
\begin{eqnarray*}
&&
\PPnmuK\QQnaB\bigl\{\|\BBhat_n-\BB\|^2 >
C_\eps n^{(1-2\beta
)/(\alpha+ 2\beta)} \bigr\}
\\
&&\qquad\le\PP_{n,\mu,K}{\Cal X}_n^c +\PP_{n,\mu,K}{
\Cal X}_n \bigl(\|\QQnaB-\QQnaBN\|_{\mathrm{TV}} + \QQnaBN
\yy_{m,\eps}^c \bigr)
\\
&&\qquad\le o_\ff(1) +2\eps+ e^{c_n} \biggl(\sum
_{i\le
n}\PP_{n,\mu,K}|\lambda_i-
\lambda_{i,N}|^2 \biggr)^{1/2}.
\end{eqnarray*}
By
construction,
\[
\lambda_i-\lambda_{i,N} = \sum
_{k> N}\zik b_k
\]
with the $\zik$'s independent and $\zik\sim N(0,\theta_k)$. Thus
\[
\sum_{i\le n}\PP_{n,\mu,K}|\lambda_i-
\lamiN|^2 \le n\sum_{k> N}
\theta_kb_k^2 = O_\ff
\bigl(nN^{1-\a-2\b}\bigr)= o_\ff\bigl(e^{-2c_n}\bigr),
\]
because $\zeta>(\a+2\b-1)^{-1}$. That is, we have an estimator that
achieves the $O_\ff(n^{(1-2\beta)/(\alpha+ 2\beta)} )$
minimax rate.

\subsection{\texorpdfstring{Proof of Theorem \protect\ref{thmeupper} with unknown Gaussian distribution}
{Proof of Theorem 1 with unknown Gaussian distribution}}\label{unknown}
Let $\BBhat_n$ be the two-stage
estimator defined in (\ref{bhat}) with cutoff points $m$ and $N$
defined in (\ref{m}) and (\ref{N}), respectively. In this section, we
show that $\BBhat_n$ achieves the asymptotic rates of convergence
stated in Theorem~\ref{thmeupper}. 

As in Section~\ref{known}, most of the analysis will be conditional
on the $\XX_i$'s lying in a set with high probability on which the
various estimators and other random quantities are well behaved. In
fact, we choose the high probability set as $\txx_{\eps, n}$ that is
defined in Lemma~\ref{txxn}. The set $\txx_{\eps,n}$ is an analogy
to ${\Cal X}_n$ in Section~\ref{known}.

As in Section~\ref{known}, the component of $\BB$ orthogonal to
$\SPAN\{\tphi_1,\ldots,\tphi_m\}$ causes no trouble because
\[
\|\BBhat-\BB\|^2 = \sum_{1\leq k\leq m}(
\bhat_k-\widetilde b_k)^2 +\bigl\|
\tH_m^\perp\BB\bigr\|^2
\]
and, by Lemma~\ref{txxn} part~\ref{txxnHm},
\[
\bigl\|\tH_m^\perp\BB\bigr\|^2 \le2\bigl\|H_m^\perp
\BB\bigr\|^2 +2\bigl\|(\tH_m-H_m)\BB\bigr\|^2 =
O_\ff\bigl(n^{(1-2\beta)/(\alpha+
2\beta
)} \bigr) \qquad\mbox{on $\txx_{\eps,n}$}.
\]
To handle 
$\sum_{1\leq k\leq m}(\bhat_k-\widetilde b_k)^2$, we invoke
Lemma~\ref{AnBn} for $\XX_i$'s lying in $\txx_{\eps,n}$, with
$\eta_i$, $A_n$, $B_n$, $D$ and $\QQ$ (notations) in Lemma
\ref{AnBn} replaced by $\teta_i$, $\tA_n$, $\tB_n$, $D$ and
$\tQQnaBN$, respectively, where
\[
\tQQnaBN:= \bigotimes_{i\le n} Q_{\tlamiN}.
\]
And, it gives a high probability set $\tyy_{m,\eps}$ with
$\tQQnaBN\tyy_{m,\eps}^c<2\eps$ on
which 
\[
\sum_{1\leq k\leq m}(\bhat_k-\widetilde
b_k)^2 = O_\ff\bigl(n^{(1-2\beta
)/(\alpha+ 2\beta)} \bigr).
\]
The conditions of
Lemma~\ref{AnBn} are satisfied on $\txx_{\eps,n}$ when $n$ is large,
because of Lemma~\ref{txxn} part~\ref{txxnteta} and
\[
\|\tA_n-\tB_n\|_2 \le\|\tA_n-
\tS A_n \tS\|_2+\|\tS A_n\tS-\tS
B_n\tS\|_2 =o_\ff(1),
\]
where the first part $\|\tA_n- \tS A_n \tS\|_2 = o_\ff(1)$ is
due to Lemma~\ref{txxn} part~\ref{txxntA}, and the second part
$\|\tS A_n\tS-\tS B_n\tS\|_2 = o_\ff(1)$ is due to Lemma
\ref{AB}.

Now, to complete the proof it suffices to show that
$\|\QQnaBN-\tQQnaBN\|_{\mathrm{TV}}$ tends to zero. First note that
\begin{eqnarray*}
\tlam_{i,N}-\lambda_{i,N} &=& a +\langle\BB,\XXbar\rangle+
\langle\tH_N\BB,\ZZ_i-\ZZbar\rangle-a-\langle\BB,\mu
\rangle- \langle H_N\BB,\ZZ_i\rangle
\\
&=& \bigl\langle\tH_N^\perp\BB,\ZZbar\bigr\rangle-\bigl
\langle H_N^\perp\BB,\ZZbar\bigr\rangle+\bigl\langle
H_N^\perp\BB,\ZZbar\bigr\rangle+ \langle\tH_N
\BB-H_N\BB,\ZZ_i\rangle,
\end{eqnarray*}
which implies that,
on $\txx_{\eps,n}$,
%
\begin{eqnarray}
\label{tlamlam}\quad
|\tlam_{i,N}-\lambda_{i,N}|^2 &\le&
2\bigl|\bigl\langle H_N^\perp\BB,\ZZbar\bigr\rangle\bigr|^2
+ 2\|\tH_N\BB-H_N\BB\|^2 \bigl(\|
\ZZ_i\| + \|\ZZbar\| \bigr)^2
\nonumber
\\
&\le& O_\ff\bigl(N^{1-2\b}\bigr)C_\eps^2
n^{-1} + O_\ff\bigl(n^{-1-\nu}\bigr)C_\eps^2
\bigl(n^{-1/2}+\rootlogn\bigr)^2
\\
&=& O_\ff\bigl(n^{-1-\nu'}\bigr)\qquad\mbox{for some $0<
\nu'<\nu$.}
\nonumber
\end{eqnarray}
Now we can argue as in step 2 of the proof for the case
of known $K$: on $\txx_{\eps,n}$,
\begin{eqnarray*}
\|\tQQnaBN-\QQnaBN\|_{\mathrm{TV}}^2 &\le& \sum
_{i\le
n}\hellinger^2 (Q_{\tlamiN},Q_{\lamiN}
)
\\
&\le& \exp\bigl(o_\ff(\log n) \bigr)\sum
_{i\le n}|\tlamiN-\lamiN|^2 \\
&=& o_\ff(1).
\end{eqnarray*}
Finish the argument as in Section~\ref{known},
by splitting into contributions from $\txx_{\varepsilon,n}^c$ and
${\txx_{\varepsilon,n}\cap\tyy_{m,\eps}^c}$ and
${\txx_{\varepsilon,n}\cap\tyy_{m,\eps}}$.

\section{\texorpdfstring{Proof of Theorem \protect\ref{thmelower}}{Proof of Theorem 2}}\label{lowerproof}
We apply a slight variation on Assouad's lemma---combining ideas from
\cite{Yu97} and from van~der Vaart [(\citeyear{vanderVaart98}), Section
24.3]---to establish the minimax lower bound result in Theorem
\ref{thmelower}.

We consider behavior only for $\mu=0$, $a=0$ and a fixed $K$ with
spectral decomposition $\sum_{j\in\NN}\theta_j \phi_j\otimes
\phi_j$
satisfying assumption {({K})}. For simplicity we abbreviate
$\PP_{n,0,K}$ to $\PP$. Let $J=\{m+1,m+2,\ldots,2m\}$ and
$\Gam=\{0,1\}^J:= \{\gamma=
(\gamma_{m+1},\ldots, \gamma_{2m}) | \gamma_j =0 \mbox{ or } \gamma_j
=1 \}$.
Let $\b_j=R j^{-\b}$. For each $\gamma$ in $\Gam$
define $\BB_\gamma= \eps\sum_{j\in J}\gamma_j\b_j \phi_j$, for a
small $\eps>0$ to be specified, and write $\QQ_\gamma$ for the product
measure $\bigotimes_{i\le n}Q_{\lambda_i(\gamma)}$ with
$
\lambda_i(\gamma)=\langle\BB_\gamma,\ZZ_i\rangle= \eps\sum_{j\in
J}\gamma_j\b_j
z_{i,j}.
$

For each $j$ let $\Gam_j=\{\gamma\in\Gam\dvtx  \gamma_j=1\}$ and let
$\psi_j$
be the bijection on $\Gam$ that flips the $j$th coordinate but
leaves all other coordinates unchanged. Let $\pi$ be the uniform
distribution on $\Gam$, that is, $\pi_\gamma=2^{-m}$ for
each $\gamma$.

For each estimator $\BBhat=\sum_{j\in\NN}\bhat_j\phi_j$ we have $
\|\BB_\gamma-\BBhat\|^2 \ge\sum_{j\in
J} (\eps\gamma_j\b_j-\bhat_j )^2 $, and so
%
\begin{eqnarray}\label{assouad}\quad
\sup_{\ff}\PP_{n,f}\|\BB-\BBhat\|^2 &\ge& \sum
_{\gamma\in
\Gam} \pi_\gamma\sum
_{j\in J}\PP\QQ_\gamma(\eps\gamma_j
\b_j-\bhat_j )^2
\nonumber
\\
&=& 2^{-m}\sum_{j\in J}\sum
_{\gamma\in\Gam_j}\PP\bigl(\QQ_\gamma(\eps\b_j-
\bhat_j)^2 +\QQ_{\psi_j(\gamma)}(0-\bhat_j)^2
\bigr)
\\
&\ge& 2^{-m}\sum_{j\in J}\sum
_{\gamma\in\Gam_j}\frac14(\eps\b_j)^2 \PP\|
\QQ_\gamma\wedge\QQ_{\psi_j(\gamma)}\|,\nonumber
\end{eqnarray}
where the first lower bound is due to the fact that the
supremum over $\ff$ is not less than the average over a subset of
$\ff$, and the last lower bound comes from the fact that
\[
(\eps\b_j-\bhat_j)^2 + (0-
\bhat_j)^2\ge\tfrac14(\eps\b_j)^2
\qquad\mbox{for all $\bhat_j$}.
\]
We assert that, if $\eps$ is chosen appropriately,
%
\begin{equation}\label{affinity}\quad
\min_{j,\gamma} \PP\|\QQ_\gamma\wedge\QQ_{\psi_j(\gamma)}\|
\qquad\mbox{stays bounded away from zero as $n\to\infty$},
\end{equation}
which will ensure that the lower bound in (\ref{assouad}) is
eventually larger than a constant multiple of $ \sum_{j\in
J}\b_j^2 \ge c n^{(1-2\beta)/(\alpha+ 2\beta)} $ for some constant
$c>0$. The inequality in
Theorem~\ref{thmelower} will then follow.

To prove (\ref{affinity}), consider a $\gamma$ in $\Gam$ and the
corresponding $\gamma'=\psi_j(\gamma)$. By virtue of the inequality
\[
\|\QQ_\gamma\wedge\QQ_{\gamma'}\| = 1-\|\QQ_\gamma-
\QQ_{\gamma
'}\|_{\mathrm{TV}} \ge1 - \biggl(2\wedge\sum
_{i\le
n}\hellinger^2(Q_{\lambda_i(\gamma)},Q_{\lambda_i(\gamma')})
\biggr)^{1/2},
\]
it is enough to show that
%
\begin{equation}\label{pph2}
\limsup_{n\to\infty}\max_{j,\gamma}\PP\biggl(2\wedge\sum
_{i\le
n}\hellinger^2(Q_{\lambda_i(\gamma)},Q_{\lambda_i(\gamma')})
\biggr) <1.
\end{equation}
Define ${\Cal X}_n = \{\max_{i\le n}\|\ZZ_i\|^2
\le C_0\log n\}$. Based on Lemma~\ref{gaussian}, we know that
$\PP{\Cal X}_n^c=o(1)$ with the constant $C_0$ large enough. On ${\Cal X}_n$
we have
\[
\bigl|\lambda_i(\gamma)\bigr|^2 \le\sum
_{j\in J}\b_j^2 \|Z_i
\|^2 = O \bigl(n^{(1-2\beta)/(\alpha+ 2\beta)}\log n \bigr) =o(1),
\]
and, by inequality in Lemma~\ref{expfactslemma}, there exits a
universal constant $C >0$ such that
\[
\hellinger^2(Q_{\lambda_i(\gamma)},Q_{\lambda_i(\gamma')}) \le C \bigl|
\lambda_i(\gamma)-\lambda_i\bigl(\gamma'
\bigr)\bigr|^2 \le C\eps^2\b_j^2
z_{i,j}^2.
\]
We deduce that
\begin{eqnarray*}
\PP\biggl( 2\wedge\sum_{i\le n} h^2(Q_{\lambda_i(\gamma)},Q_{\lambda
_i(\gamma')})
\biggr) &\le& 2\PP{\Cal X}_n^c+ C \sum
_{i\le n}\eps^2\b_j^2\PP{\Cal
X}_nz_{i,j}^2
\\
&\le& o(1)+ C \eps^2 n\b_j^2
\theta_j.
\end{eqnarray*}
The choice of $J$
makes $\b_j^2\theta_j\le R^2m^{-\a-2\b}\sim R^2/n$.
Assertion (\ref{pph2}) follows for any small enough $\eps$.

\section{Proofs of technical lemmas}\label{lemmaproof}
\subsection{\texorpdfstring{Proof of Lemma \protect\ref{AnBn}}{Proof of Lemma 1}}\label{proofAnBn} We
need to
first show the
following lemma. Note that $J_n = \sum_{i\leq n}\xi_i\xi_i'
\psiddot(\lambda_i)$. To avoid an excess of parentheses we write $\Np$
for $N+1$. We define $w_i:=J_n^{-1/2}\xi_i$ and $W_n= \sum_{i\le n}
w_i (y_i-\psidot(\lambda_i) )$.
Notice that $\QQ W_n=0$ and $ \var_\QQ(W_n)=\sum_{i\le n}
w_iw_i'\psiddot(\lambda_i) = I_\Np$ and
\[
\QQ|W_n|^2 =\trace\bigl(\var_\QQ(W_n)
\bigr) =\Np.
\]

%
\begin{lemma} \label{mleapprox} Suppose $0<\eps_1 \le1/2$ and
$0<\eps_2<1$ and
\[
{\max_{i\le n}}|w_i| \le\frac{\eps_1\eps_2}{2G(1)\Np}
\qquad\mbox{with $G$ as
in assumption {${(\dddot\Psi)}$}.}
\]
Then, the MLE $\ghat$ has the decomposition $\ghat=\gamma
+J_n^{-1/2}(W_n+r_n)$ with $|r_n|\le\eps_1$ on the
set $\{|W_n|\le\sqrt{N_+/\eps_2}\}$, which has $\QQ$-probability
greater than $1-\eps_2$.
\end{lemma}
\begin{pf}
The equality $\QQ|W_n|^2=\Np$ and Chebyshev's inequality give
\[
{\QQ\bigl\{|W_n|> \sqrt{\Np/\eps_2}\bigr\}\le\eps_2}.
\]
Reparametrize by defining $t=J_n^{1/2}(g-\gamma)$. The concave
function
\[
\ll_n(t):=L_n\bigl(\gamma+J_n^{-1/2}t
\bigr)-L_n(\gamma) =\sum_{i\le n}y_iw_i't
+\psi(\lambda_i)-\psi\bigl(\lambda_i+w_i't
\bigr)
\]
is maximized at $\that=J_n^{1/2}(\ghat-\gamma)$. It has derivative
\[
\lldot_n(t) = \sum_{i\le n}
w_i \bigl(y_i -\psidot\bigl(\lambda_i+w_i't
\bigr) \bigr).\vadjust{\goodbreak}
\]
For a fixed unit vector $u\in\RR^{N_+}$ and a fixed $t\in\RR^{N_+}$,
consider the real-valued function of the real variable $s$,
\[
H(s):= u'\lldot_n(st) = \sum
_{i\le n} u'w_i \bigl(y_i
-\psidot\bigl(\lambda_i+sw_i't\bigr)
\bigr),
\]
which has derivatives
\begin{eqnarray*}
\Hdot(s) &=& -\sum_{i\le n}\bigl(u'w_i
\bigr) \bigl(w_i't\bigr)\psiddot\bigl(
\lambda_i+sw_i't\bigr),
\\
\Hddot(s) &=& -\sum_{i\le n}\bigl(u'w_i
\bigr) \bigl(w_i't\bigr)^2\psidddot\bigl(
\lambda_i+sw_i't\bigr).
\end{eqnarray*}
Notice that $H(0)=u'W_n$ and $\Hdot(0)= -u'\sum_{i\le n} w_i
w_i'\psiddot(\lambda_i)t=-u't$.

Write $M_n:=\max_{i\le n}|w_i|$. By virtue of assumption {$
{(\dddot\Psi)}$},
\begin{eqnarray*}
\bigl|\Hddot(s)\bigr|&\le&\sum_{i\le n}\bigl|u'w_i\bigr|
\bigl(w_i't\bigr)^2 \psiddot(
\lambda_i)G \bigl(\bigl|sw_i't\bigr| \bigr)
\\
&\le& M_n G \bigl(M_n|st| \bigr)t'\sum
_{i\le n} w_iw_i'\psiddot
(\lambda_i)t
\\
&=& M_n G \bigl(M_n|st| \bigr)|t|^2.
\end{eqnarray*}
By Taylor expansion, for
some $0<s^*<1$,
\[
\bigl|H(1)- H(0)-\Hdot(0)\bigr| \le\tfrac12 \bigl|\Hddot\bigl(s^*\bigr)\bigr| \le\tfrac12
M_n G \bigl(M_n|t| \bigr)|t|^2.
\]
That is,
%
\begin{equation}\label{ull}
\bigl\llvert u' \bigl(\strut\lldot_n(t)-W_n+t
\bigr)\bigr\rrvert\le\tfrac12 M_n G \bigl(M_n|t|
\bigr)|t|^2.
\end{equation}
Approximation (\ref{ull}) will control the behavior of $\tll(s)
:=\ll_n(W_n+su)$, a concave function of the real argument $s$, for
each unit vector $u$. By concavity, the derivative $\tlldot(s)$ is a
decreasing function of $s$.
Let us decompose $\tlldot(s)$ in the following way:
\[
\tlldot(s) = u'\lldot_n(W_n+su) = -s +R(s),
\]
where
\[
\bigl|R(s)\bigr|\le\tfrac12 M_n G \bigl(M_n|W_n+su|
\bigr)|W_n+su|^2.
\]
On the set $\{|W_n|\le\sqrt{\Np/\eps_2}\}$ we have
\[
|W_n\pm\eps_1 u|\le\sqrt{\Np/\eps_2} +
\eps_1.
\]
Thus
\[
M_n |W_n\pm\eps_1 u|\le\frac{\eps_1\eps_2}{2G(1)\Np} (
\sqrt{\Np/\eps_2} +\eps_1 )<1,
\]
implying
\begin{eqnarray*}
\bigl|R(\pm\eps_1)\bigr| &\le& \frac12 M_n G(1)|W_n\pm
\eps_1 u|^2 \le\frac{\eps_1\eps_2}{G(1)\Np} \bigl(\Np/
\eps_2 +\eps_1^2 \bigr)
\\
&\le& \eps_1 \bigl(1+\eps_1^2
\eps_2/\Np\bigr) <\frac58\eps_1.
\end{eqnarray*}
Deduce that
\[
\tlldot(\eps_1) = -\eps_1+R(\eps_1)\le-
\tfrac38\eps_1 \quad\mbox{and}\quad \tlldot(-\eps_1) =
\eps_1+R(-\eps_1)\ge\tfrac38\eps_1.
\]
The
concave function $s\mapsto\ll_n(W_n+su)$ must achieve its maximum
for some $s$ in the interval $[-\eps_1,\eps_1]$, for each unit
vector $u$. It follows that $|\that-W_n|\le\eps_1$.
\end{pf}
First we establish a bound on the spectral distance between $A_n^{-1}$
and $B_n^{-1}$. Define $H=B_n^{-1}A_n-I$.
Then$\|H\|_2\le\|B_n^{-1}\|_2\|A_n-B_n\|_2\le1/2$, which
justifies the expansion
\[
\bigl\|A_n^{-1}-B_n^{-1}\bigr\|_2
= \bigl\| \bigl((I+H)^{-1}-I \bigr)B_n^{-1}
\bigr\|_2 \le\sum_{j\ge1}\|H
\|_2^k\bigl\|B_n^{-1}\bigr\|_2 \le
\bigl\|B_n^{-1}\bigr\|_2.
\]
As a consequence, $\|A_n^{-1}\|_2 \le2\|B_n^{-1}\|_2$.

Choose $\eps_1=1/2$ and $\eps_2=\eps$ in Lemma~\ref{mleapprox}. The
bound on $\max_{i\le n}|\eta_i|$ gives the bound on $\max_{i\le
n}|w_i|$ needed by the lemma
\[
n|w_i|^2 = \eta_i'D(J_n/n)^{-1}D
\eta_i = \eta_i'A_n^{-1}
\eta_i \le\bigl\|A_n^{-1}\bigr\|_2|
\eta_i|^2.
\]
As shown in Lemma~\ref{mleapprox}, the MLE $\hat{g}$ can be
decomposed as
\[
\hat{g} = \gamma+ J_n^{-1/2}(W_n +
r_n).
\]
Define $K_j:=J_n^{-1/2}\kappa_j$, so that $
|\kappa_j'(\ghat-\gamma)|^2 \le2(K_j'W_n)^2 + 2(K_j'r_n)^2 $. By
Cauchy--Schwarz,
\[
\sum_{j} \bigl(K_j'r_n
\bigr)^2 \le\sum_{j}|K_j|^2
|r_n|^2 = U_\kappa|r_n|^2,
\]
where
\[
U_\kappa:=\sum_{j}\kappa_j'J_n^{-1}
\kappa_j = \sum_{j}n^{-1}
\bigl(D^{-1}\kappa_j\bigr)'A_n^{-1}D^{-1}
\kappa_j 
\le2n^{-1}\bigl\|B_n^{-1}
\bigr\|_2 \sum_{j}\bigl|D^{-1}
\kappa_j\bigr|^2.
\]
For the\vspace*{2pt} contribution $V_\kappa:= \sum_{j}|K_j'W_n|^2$, the
Cauchy--Schwarz bound is too crude. Instead, notice that $\QQ
V_\kappa=U_\kappa$, which ensures that the complement of the set
\[
\yy_{\kappa,\eps}:= \bigl\{|W_n|\le\sqrt{\Np/\eps}\bigr\} \cap
\{V_\kappa\le U_\kappa/\eps\}
\]
has $\QQ$ probability less that $2\eps$. On the
set $\yy_{\kappa,\eps}$,
\[
\sum_{0\le j\le N}\bigl|\kappa_j'(
\ghat-\gamma)\bigr|^2 \le2V_\kappa+2U_\kappa|r_n|^2
\le3U_\kappa/\eps.
\]
The asserted bound follows.

\subsection{\texorpdfstring{Proof of Lemma \protect\ref{AB}}{Proof of Lemma 2}}\label{proofAB}
Throughout this subsection, abbreviate $\PP_{n,\mu,K}$ to $\PP$.
The matrix $A_n$ is an average of $n$ independent random matrices
each of which is distributed like ${\Cal N}{\Cal N}'\psiddot(\gamma
'D{\Cal N})$,
where ${\Cal N}=({\Cal N}_0,{\Cal N}_1,\ldots,{\Cal N}_N)'$ with ${\Cal
N}_0\equiv1$, and the
other ${\Cal N}_j$'s are independent $N(0,1)$'s. Moreover, by rotational
invariance of the spherical normal, we may assume with no loss of
generality that $\gamma'D{\Cal N}= \abar+\kappa{\Cal N}_1$, where
\[
\kappa^2 = \sum_{k=1}^N
D_k^2 b_k^2 =
O_\ff(1).
\]
Thus
\[
B_n =\PP{\Cal N} {\Cal N}'\psiddot(\abar+\kappa{\Cal
N}_1) = \diag(F, r_0 I_{N-1}),
\]
where
\[
r_j:= \PP{\Cal N}_1^j\psiddot(\abar+\kappa{
\Cal N}_1) \quad\mbox{and}\quad F = %
\left[\matrix{r_0&r_1
\cr
r_1&r_2} \right].
\]
The block diagonal form of $B_n$ simplifies calculation of spectral
norms,
\begin{eqnarray*}
\bigl\|B_n^{-1}\bigr\|_2 &=& \bigl\|\diag\bigl(F^{-1},r_0^{-1}I_{N-1}
\bigr)\bigr\|_2
\\
&\le&\max\bigl(\bigl\|F^{-1} \bigr\|_2, \bigl\|r_0^{-1}I_{N-1}
\bigr\|_2 \bigr) \le\max\biggl(\frac{r_0+r_2}{r_0r_2-r_1^2}, r_0^{-1}
\biggr).
\end{eqnarray*}
Assumption {$({\ddot\Psi})$} ensures that both $r_0$ and $r_2$ are
$O_\ff(1)$.

Continuity and strict positivity of $\psiddot$, together with
$\max(|\abar|,\kappa)=O_\ff(1)$, ensure that $ c_0:=
\inf_{\abar,\kappa}\inf_{|x|\le1}\psiddot(\abar+\kappa x)>0
$. Thus
\[
\sqrt{2\pi} r_0 \ge c_0\int_{-1}^{+1}
e^{-x^2/2}\,dx >0.
\]
Similarly,
\begin{eqnarray*}
\sqrt{2\pi}\bigl(r_0r_2-r_1^2
\bigr) &=& \sqrt{2\pi}r_0\PP\psiddot(\abar+\kappa{\Cal
N}_1) ({\Cal N}_1-r_1/r_0)^2
\\
&\ge& c_0r_0\int_{-1}^{+1}
(x-r_1/r_0)^2e^{-x^2/2}\,dx \\
&\ge&
c_0r_0\int_{-1}^{+1}
x^2e^{-x^2/2}\,dx.
\end{eqnarray*}
It follows that
$\|B_n^{-1}\|_2=O_\ff(1)$.

The random matrix $A_n-B_n$ is an average of $n$ independent random
matrices, each distributed like ${\Cal N}{\Cal N}'\psiddot(\abar
+\kappa{\Cal N}_1)$
minus its expected value. Thus
\[
\PP\|A_n-B_n\|_2^2\le\PP
\|A_n-B_n\|_{\mathrm{F}}^2 =
n^{-1}\sum_{0\le
j,k \le N}\var\bigl({\Cal
N}_j{\Cal N}_k \psiddot(\abar+\kappa{\Cal
N}_1) \bigr),
\]
where $\|\cdot\|_{\mathrm{F}}$ is the Frobenius norm.
Assumption {$({\ddot\Psi})$} ensures
that each summand is $O_\ff(1)$, which
leaves us with a $O_\ff(N^2/n) =o_\ff(1)$ upper bound.

\subsection{\texorpdfstring{Proof of Lemma \protect\ref{expfactslemma}}{Proof of Lemma 3}}\label{hellinger}
Let us temporarily write $\lambda'$ for $\lambda+\delta$ and write
$\barlam$
for ${(\lambda+\lambda')/2=\lambda+\delta/2}$,
\begin{eqnarray*}
1-\frac12\hellinger^2(Q_\lambda,Q_{\lambda'}) &=& \int
\sqrt{f_\lambda(y)f_{\lambda
'}(y)}= \int\exp\biggl(\barlam y -
\frac12\psi(\lambda)-\frac12\psi\bigl(\lambda'\bigr) \biggr)
\\
&=& \exp\biggl(\psi(\barlam) -\frac12\psi(\lambda)-\frac12\psi
\bigl(\lambda'\bigr) \biggr)\\
&\ge& 1 + \psi(\barlam) - \frac12\psi(
\lambda) -\frac12\psi\bigl(\lambda'\bigr).
\end{eqnarray*}
That is,
\[
\hellinger^2(Q_\lambda,Q_{\lambda'})\le\psi(\lambda)+
\psi(\lambda+\delta) - 2\psi(\lambda+\delta/2).
\]
By Taylor expansion in $\delta$ around $0$, the right-hand side is
less than
\[
\tfrac14\delta^2
\psiddot(\lambda) +\tfrac16\delta^3 \bigl(\psidddot\bigl(\lambda+
\delta^*\bigr)-\tfrac14\psidddot\bigl(\lambda+ \delta^*/2\bigr) \bigr),
\]
where $0<|\delta^*|<|\delta|$. Invoke inequality {${(\dddot
\Psi)}$} twice
to bound the coefficient of $\delta^3/6$ in absolute value by
\[
\psiddot(\lambda) \bigl(\strut G\bigl(|\delta|\bigr) +\tfrac14 G\bigl(|\delta|/2\bigr)
\bigr) \le
\tfrac54\psiddot(\lambda)G\bigl(|\delta|\bigr).
\]
The stated bound simplifies some unimportant constants.

\subsection{\texorpdfstring{Proof of Lemma \protect\ref{gaussian}}{Proof of Lemma 4}}\label{proofgaussian}
Without loss of generality, let us suppose $T=1$. For $s=1/4$, note that
\[
\PP\bigl[ \exp(sW_i) \bigr] =\prod_{k\in\NN}(1-2s
\tau_{i,k})^{-1/2} 
\le\exp\biggl(\sum
_{k\in\NN}s\tau_{i,k} \biggr)\le e^{1/4}
\]
by
virtue of the inequality $-\log(1-t)\le2t$ for $|t|\le1/2$. With
the same $s$, it then follows that
\begin{eqnarray*}
\PP\Bigl\{\max_{i\le n} W_i > 4(\log n + x)\Bigr\}&\le&\exp
\bigl(-4s(\log n + x) \bigr)\PP\Bigl[\exp\Bigl(\max_{i\le n}
sW_i \Bigr) \Bigr]
\\
&\le& e^{-x}\frac1n\sum_{i\le n}\PP\bigl[
\exp(sW_i) \bigr].
\end{eqnarray*}
The $2$ is just a clean upper bound for $e^{1/4}$.

\subsection{\texorpdfstring{Proof of Lemma \protect\ref{txxn}}{Proof of Lemma 5}}\label{prooftxxn}

We first show some preliminary lemmas in Section~\ref{txxnPrelim}.
Those preliminary results are used in the main proofs throughout
Sections~\ref{txxnXX} to~\ref{txxnSix}. For notational simplicity, we
write
$\sum_j^*$ for $\sum_{j\neq k}$. 

\subsubsection{Preliminary lemmas}\label{txxnPrelim}

Remember that $\theta_j$'s are the eigenvalues of $K$ as defined in
Definition~\ref{ffdef}. Many of the
inequalities in the proof of Lemma~\ref{txxn} involve sums of
functions of the $\theta_j$'s. The following result will save us a lot
of repetition.\vadjust{\goodbreak}
%
\begin{lemma} \label{weights} 
\textup{(i)} For each $r\ge1$ there is a constant $C_r=C_r(\ff)$ for
which
\[
\kappa_k(r,\gamma):= \sum_{j\in\NN}\{j\ne
k\}\frac{j^{-\gamma}}{|\theta_j-\theta_k|^r} \le%
\cases{C_r
\bigl(1+k^{r(1+\a)-\gamma} \bigr), &\quad if $r>1$,
\vspace*{2pt}\cr
C_1
\bigl(1+k^{1+\a-\gamma}\log k \bigr), &\quad if $r=1$. } %
\]

\textup{(ii)} For each $p$,
\[
\sum_{k\le p}\sum_{j>p}
\frac{ k^{-\a-2\b}j^{-\a} }{
|\theta_k-\theta_j|^2} = O_\ff\bigl(p^{1-\a}\bigr).
\]
\end{lemma}
\begin{pf} For (i), argue in the same way as
Hall and Horowitz [(\citeyear{HallHorowitz2007}), pa\-ge~85], using the lower bounds
\[
|\theta_j-\theta_k| \ge%
\cases{
c_\a j^{-\a}, &\quad if $j<k/2$,
\cr
c_\a|j-k|
k^{-\a-1}, &\quad if $k/2\le j\le2k$,
\cr
c_\a k^{-\a}, &\quad
if $j>2k$, } %
\]
where $c_\a$ is a positive constant.

For (ii), split the range of summation into two subsets: $\{(k,j)\dvtx
j> \break\max(p,2k)\}$ and $\{(k,j)\dvtx  p/2<k\le p<j\le2k\}$. The first
subset contributes at most
\[
\sum_{k\le p} k^{-\a-2\b}\sum
_{j> \max(p,2k)}j^{-\a}\bigl(c_\a k^{-\a}
\bigr)^{-2} = O_\ff\bigl(p^{1-\a}\bigr),
\]
because $\a-2\b<-3$. The second subset contributes at most
\[
\sum_{p/2<k\le
p}k^{-\a-2\b}c_\a^{-2}k^{2\a+2}
\sum_{j>p}j^{-\a}(j-k)^{-2}
=O_\ff\bigl(p^{2+\a-2\b}p^{1-\a} \bigr),
\]
which is of order $o_\ff(p^{-\a})$.
\end{pf}

Remember that $\zij=\langle\ZZ_i,\phi_j\rangle$ and the standardized
variables $\etaij=\zij/\break\sqrt{\theta_j}$ are independent $N(0,1)$'s.
Define $\eta_{\cdot j}=n^{-1}\sum_{i\le n}\etaij$ and
\[
\scov_{j,k}:=(n-1)^{-1} \sum_{i\le
n}
(\etaij-\eta_{\cdot j} ) (\etaik-\eta_{\cdot
k} ),
\]
the $(j,k)$-element of a sample covariance matrix of
i.i.d. $N(0,I_N)$ random vectors. We further define
\[
\Lam_k:= \sum_{j\in\NN}
\Lam_{k,j}\phi_j \qquad\mbox{with } \Lam_{k,j}:=
\cases{ \sqrt{\theta_j\theta_k}
\scov_{j,k}/(\theta_k-\theta_j), &\quad if $j\ne k$,
\cr
0, &\quad if $j=k$. } %
\]

In fact, most of the inequalities that we need for proving Lemma
\ref{txxn} come from simple moment bounds (Lemma~\ref{Sjkfacts})
for the sample covariances $\scov_{j,k}$ and the derived bounds
(Lemma~\ref{Lambounds}) for the $\Lam_k$'s. The distribution
of $\scov_{j,k}$ does not depend on the parameters of our model.
By the rotation of axes we can rewrite
$(n-1)\scov_{j,k}$ as $U_j'U_k$, where $U_1,U_2,\ldots$ are
independent $N(0,I_{n-1})$ random vectors. This representation gives
some useful equalities and bounds.\vadjust{\goodbreak}

\begin{lemma}\label{Sjkfacts} 
Uniformly over distinct $j,k,\ell$:
\begin{longlist}
\item
$\PP\scov_{j,j}= 1$ and $\PP(\scov_{j,j}-1 )^2 =
2(n-1)^{-1};$

\item
$\PP\scov_{j,k} = \PP\scov_{j,k} \scov_{j,\ell}=0; $

\item
$\PP\scov_{j,k}^2 =O(n^{-1}).$
\end{longlist}
\end{lemma}
\begin{pf} Assertion (i) is classical because $|U_j|^2\sim\chi_{n-1}^2$.
For assertion (ii) use $ \PP(U_1'U_2| U_2) =0 $ and
\[
\PP\bigl(U_1'U_2U_2'U_3
| U_2\bigr) = \trace\bigl(U_2U_2'
\PP\bigl(U_3U_1'\bigr) \bigr)=0.
\]
For (iii) use $\PP(U_1U_1')=I_{n-1}$ and
\[
\PP\bigl(U_1'U_2U_2'U_1
| U_2\bigr) = \trace\bigl(U_2U_2'
\PP\bigl(U_1U_1'\bigr) \bigr)= \trace
\bigl(U_2U_2'\bigr) =|U_2|^2.
\]
\upqed\end{pf}
%
\begin{lemma} \label{Lambounds} 
Uniformly over distinct $j,k,\ell$:
{\renewcommand\thelonglist{(\roman{longlist})}
\renewcommand\labellonglist{\thelonglist}
\begin{longlist}
\item\label{Lambounds0}
$\PP\Lam_{k,j} = \PP\Lam_{k,j} \Lam_{k,\ell}=0$;

\item
$\PP\Lam_{k,j}^2 = O_\ff(n^{-1} k^{-\a}j^{-\a}
(\theta_k-\theta_j)^{-2} )$;


\item\label{LamboundsLam2}
$\PP\|\Lam_k\|^2 =O_\ff(n^{-1}k^2)$.
%
\end{longlist}}
\end{lemma}
\begin{pf} Assertions (i) and (ii) follow from assertions (ii)
and (iii) of Lem\-ma~\ref{Sjkfacts}. For (iii), note that
\[
\PP\|\Lam_k\|^2 = \sum_j^*
\PP\Lam_{j,k}^2 = O_\ff\bigl(n^{-1}k^{-\a}
\bigr)\kappa_k(2,\a)
\]
and $\kappa_k(2,\a) = O_\ff(k^{2+\a})$ from Lemma~\ref{weights}.
\end{pf}

The following two lemmas related to perturbation theory for
self-adjoint compact operators [cf., e.g., \cite
{Bosq2000,BirmanSolomjak87,Kato95}] are crucial in the development of Lemma
\ref{txxn}. They are special cases of Lemmas 2 
and 4 
in the supplemental article [\cite{DouPollardZhou2012Supplemental}] under
the general
perturbation-theoretic framework. For Lemma~\ref{eigenvec}, similar
results were established by other authors; see, for example, Hall and
Hosseini-Nasab (\citeyear{HallHosseini2006}), equation 2.8, and
\cite{CaiHall2006}, Section 5.6. Lemma~\ref{eigenspace} extends the perturbation
result for eigenprojections, obtained by Tyler [(\citeyear{Tyler81}), Lemma 4.1],
from the matrix case to the general operator case.

Define
\[
\varepsilon_k:= \min\bigl\{|\theta_j - \theta_k|\dvtx j
\neq k\bigr\}
\]
and
\[
f_k:= \sigma_k \widetilde{\phi}_k -
\phi_k\qquad\mbox{for all } k.
\]

\begin{lemma}\label{eigenvec} 
If we have $\varepsilon_k>5\|\Delta\|$, then it follows that
\[
\|f_k\|\leq3\|\Lam_k\|.\vadjust{\goodbreak}
\]
\end{lemma}

Define $H_J = \SPAN\{\phi_j\dvtx j\in J\}$ and $\tH_J =
\SPAN\{\tphi_j\dvtx j\in J\}$ for $J\subseteq\NN$.
%
\begin{lemma}\label{eigenspace} 
If we have
$\min_{k\in J}\eps_k>5\|\Delta\|$, then it follows that
\[
(\tH_J-H_J)\BB= \sum_{j\in
J}
\sum_{k\in J^c}\phi_j b_k (
\Lam_{j,k} + \Lam_{k,j}) + e,
\]
where
$\|e\|^2$ is bounded by a universal constant times $R_1 +
\|\Delta\|^2 R_2$ with
\begin{eqnarray*}
R_1 &=& \biggl(\sum_{k\in
J}\|
\Lam_k\|^2 \biggr) \sum_{k\in
J}
\Biggl(\sum_{j}^*\Lam_{k,j}b_j
\Biggr)^2,
\\
R_2&=&\sum_{k\in
J}\|\Lam_k
\|^2 \Biggl(\sum_{j}^*\frac{|b_j|}{|\theta_k-\theta_j|}
\Biggr)^2 + \Biggl(\sum_{k\in
J}\|
\Lam_k\||b_k|\sum_{j}^*
\frac{1}{|\theta_k-\theta_j|} \Biggr)^2
\\
&&{} +\sum_{k\in J}\|\Lam_k
\|^2|b_k|^2k^{2+2\alpha}.
\end{eqnarray*}
\end{lemma}

\subsubsection{\texorpdfstring{A high probability set $\txx_{\eps, n}$}
{A high probability set X epsilon, n}}\label{txxnXX}
To prove Lemma~\ref{txxn} we define $\txx_{\eps,n}$ as an
intersection of sets chosen to make the six assertions of the lemma
hold,
\[
\txx_{\eps,n}:= \txx_{\Del,n}\cap\txx_{\ZZ,n}\cap
\txx_{\eta,n}\cap\txx_{A,n} \cap\txx_{\Lam,n},
\]
where the complement of each of the five sets appearing on the
right-hand side has probability less than $\eps/5$. More
specifically, for a large enough constant $C_\eps$, we define
\begin{eqnarray*}
\txx_{\Del,n} &=& \bigl\{\|\Del\|\le C_\eps n^{-1/2}
\bigr\},
\\
\txx_{\ZZ,n} &=& \Bigl\{{\max_{i\le n}}\|\ZZ_i
\|^2 \le C_\eps\log n\mbox{ and } \|\ZZbar\|\le
C_\eps n^{-1/2}\Bigr\},
\\
\txx_{\eta,n} &=& \biggl\{{\max_{i\le n}}|\eta_i|^2
\le C_\eps N\log n \mbox{ and } \biggl\|\sum_{i\leq n}
\eta_i\eta_i'\biggr\|_2 \leq
C_\varepsilon n\biggr\}, 
\\
\txx_{A,n} &=& \biggl\{\biggl\|\sum_{i\le n}
\teta_i\teta_i' \biggr\|_2 \le
C_\eps n\biggr\}.
\end{eqnarray*}
The set of $\txx_{\Lam,n}$ is defined in a slightly more
complicated way. It is defined by requiring various functions of
the $\Lam_k$'s to be smaller than $C_\eps$ times their expected
values. Calculate expected values for all the terms in $R_1$ and
$R_2$ that appear in the bound of Lemma~\ref{eigenspace}.
%
\begin{eqnarray}\label{proj1}
&&
\PPnmuK\sum_{k\le p} \biggl(\sum
_{j>p}\Lam_{k,j}b_j
\biggr)^2 + \PPnmuK\sum_{j>p} \biggl(\sum
_{k\le p}\Lam_{k,j}b_k
\biggr)^2
\nonumber
\\
&&\qquad= \sum_{k\le
p}\sum_{j>p}
\PPnmuK\Lam_{k,j}^2 \bigl(b_j^2+b_k^2
\bigr) \qquad\mbox{by Lemma~\ref{Lambounds} part~\ref{Lambounds0}}
\\
&&\qquad= O_\ff\bigl(n^{-1} \bigr)\sum
_{k\le
p}\sum_{j>p}k^{-\a-2\b}j^{-\a}(
\theta_k-\theta_j)^{-2}
\nonumber
\\
&&\qquad= O_\ff\bigl(n^{-1}p^{1-\a} \bigr) \qquad\mbox{by
Lemma~\ref{weights}}\nonumber
\end{eqnarray}
and
%
\begin{eqnarray*}
\PPnmuK\sum_{k\le p}b_k^2\|
\Lam_k\|^2 k^{2+2\alpha} &=& O_\ff
\bigl(n^{-1} \bigr)\sum_{k\le p}k^{4+2\alpha-2\b}
\\
&=& O_\ff\bigl(n^{-1} \bigr) \bigl(1+p^{5+2\alpha-2\b}+
\log p \bigr)
\end{eqnarray*}
and
\[
\PPnmuK\sum_{k\le p}|b_k|\|
\Lam_k\|^2 = O_\ff\bigl(n^{-1}
\bigr)\sum_{k\in J}k^{2-\b} = O_\ff
\bigl(n^{-1} \bigr) \bigl(1+p^{3-\b}+\log p \bigr)
\]
and
\[
\PPnmuK\sum_{k\le p}\|\Lam_k
\|^2 = O_\ff\bigl(n^{-1}p^3 \bigr)
\]
and
%
\begin{eqnarray}\label{proj2}\quad
\PPnmuK\sum_{k\le
p} \Biggl(\sum
_{j}^*\Lam_{k,j}b_j
\Biggr)^2 &=& O_\ff\bigl(n^{-1} \bigr)\sum
_{k\le p}\sum_{j}^*
k^{-\a}j^{-a-2\b}(\theta_k-\theta_j)^{-2}
\nonumber\\[-8pt]\\[-8pt]
&=& O_\ff\bigl(n^{-1} \bigr) \qquad\mbox{by Lemma
\ref{weights}}\nonumber
\end{eqnarray}
and
%
\begin{equation}\label{proj3}\qquad
\PPnmuK\sum_{k\le p} \|\Lam_k
\|^2 \Biggl(\sum_j^*\frac{|b_j|}{|\theta_k-\theta_j|}
\Biggr)^2 = O_\ff\bigl(n^{-1} \bigr)
\bigl(p^3+ p^{5+2\a-2\b}\log^2 p \bigr)
\end{equation}
and by
Lemma~\ref{weights}
%
\begin{equation}\label{proj4}
\sum_{k\le
p}b_k^2 \Biggl(\sum
_j^*\frac{1}{|\theta_k-\theta_j|} \Biggr)^2
=O_\ff\bigl(1+p^{3+2\a-2\b}\log^2 p \bigr).
\end{equation}
For some constant
$C_\eps=C_\eps(\ff)$, on a set ${\Cal X}_{\Lam,n}$ with $\PPnmuK
{\Cal X}_{\Lam,n}^c <\eps/5$, each of the random quantities in the
previous set of inequalities (for both $p=m$ and $p=N$) is bounded
by $C_\eps$ times its $\PPnmuK$ expected value. By virtue\vspace*{1pt} of
Lemma~\ref{Lambounds} part~\ref{LamboundsLam2}, we may also
assume that $\|\Lam_k\|^2\le C_\eps k^2/n$ on ${\Cal
X}_{\Lam,n}$.\vspace*{1pt}

We now show that $\sup_{f\in\ff}\PPnmuK\txx_{\eps, n}^c < \eps$.
From the construction of $\txx_{\Lam, n}$ above, it follows directly
that $\PPnmuK\txx_{\Lam, n}^c < \eps/5$.

We analyze $\tK$ by rewriting it using the eigenfunctions for $K$.
Then
\[
\ZZ_i(t)-\ZZbar(t) = \sum_{j\in\NN} (
\zij-z_{\cdot j})\phi_j(t) = \sum_{j\in\NN}
\sqrt{\theta_j}(\eta_{i,j}-\eta_{\cdot j})
\phi_j(t)
\]
and
%
\begin{equation}\label{Khatrep}
\tK(s,t) = \sum_{j,k\in\NN}\tK_{j,k}
\phi_j(s)\phi_k(t) \qquad\mbox{with $\tK_{j,k} =
\sqrt{\theta_j\theta_k}\scov_{j,k}$}.
\end{equation}
Observe that
\begin{eqnarray*}
\PP\|\Del\|^2 &=& \sum_{j,k}\PPnmuK
\bigl(\tK_{j,k}-\theta_j\{j=k\} \bigr)^2 = \sum
_{j,k}\theta_j\theta_k\PP
\bigl(\scov_{j,k}-\{j=k\} \bigr)^2
\\
&\le& \sum_{j}\theta_jO_\ff
\bigl(n^{-1}\bigr) + \sum_{j,k}
\theta_j\theta_k O_\ff\bigl(n^{-2}
\bigr)= O_\ff\bigl(n^{-1}\bigr).
\end{eqnarray*}
Thus, we have $\PPnmuK
\txx_{\Delta, n}^c < \eps/5$.

The set $\txx_{A,n}$ is almost redundant in the sense that
$\txx_{\Delta, n} \subseteq\txx_{A,n}$ when $n$ and $C_{\eps}$ are
large enough. From Definition~\ref{ffdef} we know that
\[
{\min_{1\le j<j'\le N}}|\theta_j-\theta_{j'}| \ge(\a/R)
N^{-1-\a} \quad\mbox{and}\quad\min_{1\le j\le N}\theta_j \ge
R^{-1}N^{-\a}.
\]
The choice $N\asymp n^\zeta$ with $\zeta< (2+2\a)^{-1}$ ensures
that $n^{1/2}N^{-1-\a}\to\infty$. On $\txx_{\Del,n}$ the spacing
assumption used in Lemmas~\ref{eigenvec} and~\ref{eigenspace} holds
for all $n$ large enough; all the bounds from those lemmas are
available to us on $\txx_{\eps,n}$. In particular,
%
\begin{equation}\label{maxth}
{\max_{j\le
N}}|\tth_j/\theta_j-1| \le
O_\ff\bigl(N^\a\|\Del\|\bigr) =o_\ff(1),
\end{equation}
where $\tth_j$'s are defined in (\ref{tdecomp}).
Remember that
\[
\ZZ_i(t) -\ZZbar(t) = \sum_{k\in\NN}(
\tzik-{\widetilde z_{\cdot
k}})\tphi_k(t)
\]
so that
\[
\tth_k\{j=k\} = \iint\tK(s,t)\tphi_j(s)
\tphi_k(t)\,ds\,dt =(n-1)^{-1}\sum
_{i\le
n}(\tz_{i,j}-{\widetilde z_{\cdot j}}) (
\tz_{i,k}-{\widetilde z_{\cdot
k}}),
\]
which implies
$(n-1)^{-1}\sum_{i\le n} \tz_i\tz_i' = \tD^2$ and
%
\begin{equation}\label{tetateta}
(n-1)^{-1}\sum_{i\le n} \teta_i
\teta_i' = D^{-1}\tD^2D^{-1}
:= \diag(1,\tth_1/\theta_1,\ldots,\tth_N/
\theta_N).
\end{equation}
Inequality (\ref{maxth}) and equality (\ref{tetateta}) together show
that $\txx_{\Del,n}\subseteq\txx_{A,n}$ eventually if we make sure
$C_\eps>1$. Thus, $\PPnmuK\txx_{A, n}^c \leq\PPnmuK\txx_{\Delta,
n}^c < \eps/5$.

As the controls for the set defined in (\ref{xxnZZ}) and
(\ref{xxneta}), Lemma~\ref{gaussian} controls ${\max_{i\le
n}}\|\ZZ_i\|^2$ and ${\max_{i\leq n}}|\eta_i|^2$. In addition, we
know that
\begin{eqnarray*}
\PP\biggl\|n^{-1}\sum_{i\leq n} \eta_i
\eta_i' - I_{N+1}\biggr\|_2^2
&\leq& \PP\biggl\|n^{-1}\sum_{i\leq n}
\eta_i\eta_i' - I_{N+1}
\biggr\|_{\mathrm{F}}^2
\\
&=& n^{-1}\sum_{0\leq j,k\leq N}\operatorname{var} (
\eta_{i,k}\eta_{i,j} )=O_\ff\bigl(N^2/n
\bigr).
\nonumber
\end{eqnarray*}
Thus, $\PP\|n^{-1}\sum_{i\leq n} \eta_i\eta_i'\|_2 = 1 +
o_\ff
(1)$. Therefore, we have
$\PPnmuK\txx_{\eta, n}^c < \eps/5$. To control the $\ZZbar$
contribution, note that $n\|\ZZbar\|^2$ has the same distribution
as $\|\ZZ_1\|^2$, which has expected value
$\sum_{j\in\NN}\theta_j<\infty$. Thus, we have $\PPnmuK\txx_{\ZZ,
n}^c < \eps/5$.

Therefore, there exists $C_\eps>0$ such that
\[
\PPnmuK\txx_{\eps,
n}^c \leq\PPnmuK\bigl(
\txx_{\Delta, n}^c + \txx_{\ZZ, n}^c +
\txx_{\eta, n}^c + \txx_{A, n}^c +
\txx_{\Lam, n}^c \bigr)< \eps.
\]

\subsection{\texorpdfstring{Proof of the assertions on $\txx_{\eps, n}$}{Proof of the assertions on X epsilon, n}}\label{txxnSix}

The assertions~\ref{txxnDel} and~\ref{txxnZZ} hold on the set
$\txx_{\eps, n}$ as a direct consequence of the construction. From
Lemma~\ref{eigenspace}, it follows that on the set ${\Cal
X}_{\Del,n}\cap {\Cal X}_{\Lam,n}$, if $p\le N$,
%
\[
\bigl\|(\tH_p-H_p)\BB\bigr\|^2 = O_\ff
\bigl(n^{-1}p^{1-\a} \bigr).
\]
This inequality leads to the asserted conclusions in
\ref{txxnHm} and~\ref{txxnHN} when $p=m$ or $p=N$.

Now we show assertion~\ref{txxnteta} holds on the set
$\txx_{\eps, n}$. By construction, $\teta_{i1}=1$ for every $i$, and
for $j\ge2$,
\[
\sqrt{\theta_j}\teta_{i,j} =(\tz_{i,j}-{\widetilde
z_{\cdot j}}) = \langle\ZZ_i-\ZZbar,\tphi_j
\rangle.
\]
Thus, for $j\ge2$,
\[
\sig_j\teta_{i,j} = \theta_j^{-1/2}
\langle\ZZ_i-\ZZbar,\phi_j+f_j\rangle=
\eta_{i,j}+\tdel_{i,j}
\]
with $\tdel_{i,j}$ satisfying the following bound, due to Lemma \ref
{eigenvec}:
\[
|\tdel_{i,j}|^2\le\theta_j^{-1} \bigl(\|
\ZZ_i\|+\|\ZZbar\| \bigr)^2\|f_j\|^2
\le O_\ff\biggl(\frac{j^{2+\a}\log n}{n} \biggr) \qquad\mbox{on $
\txx_{\eps,n}$}.
\]
In vector form,
%
\begin{equation}\label{tetaeta}
\tS\teta_i = \eta_i + \tdel_i \qquad\mbox{with
}|\tdel_i|^2 =O_\ff\biggl(
\frac{N^{3+\a}\log n}{n} \biggr)\le o_\ff\bigl(n/N^2\bigr)
\qquad\mbox{on }\txx_{\eps,n}.\hspace*{-35pt}
\end{equation}
It
follows that
\[
{\max_{i\le n}}|\teta_i| = {\max_{i\le n}}|\tS
\teta_i| \le{\max_{i\le n}}|\eta_i|+o_\ff(
\sqrt{n}/N) = O_\ff(\sqrt{n}/N) \qquad\mbox{on $\txx_{\eps,n}$}.
\]

In the end, we show that on $\txx_{\eps, n}$ assertion
\ref{txxntA} holds. From inequality (\ref{tlamlam}) we know that
\[
\tilde{\eps}_{N}:= \max_{i\le n}|\tlam_{i,N}-
\lambda_{i,N}| = O_\ff\bigl(n^{-(1+\nu')/2}\bigr) \qquad\mbox{on
$\txx_{\eps,n}$},
\]
and from bounds (\ref{deltai}) and (\ref{lambdai}) in Section \ref
{known}, we have $ \max_{i\le
n}|\lambda_{i,N}| = O_\ff(\rootlogn) $. Assumption {${(\dddot
\Psi)}$} in
Section~\ref{model} and the mean-value theorem then give
\[
\max_{i\le n}\bigl|\psiddot(\tlam_{i,N})-\psiddot(
\lambda_{i,N})\bigr| \le\tilde{\eps}_N\psiddot(
\lambda_{i,N})G(\tilde{\eps}_N) =o_\ff(1).
\]
If we replace $\psiddot(\tlam_{i,N})$ in the definition of $\tA_n$
by $\psiddot(\lambda_{i,N})$, we make a change
\[
\Pi= (n-1)^{-1}\sum_{i\leq n}\tilde{
\eta}_i\tilde{\eta}_i'\bigl(\psiddot(
\tlam_{i,N}) - \psiddot(\lambda_{i,N})\bigr)
\]
with
$
\|\Pi\|_2 \le o_\ff(1)\|(n-1)^{-1}\sum_{i\le n} \teta_i\teta_i' \|_2,
$
which, by equality (\ref{tetateta}), is of order $o_\ff(1)$
on $\txx_{\eps,n}$.

From assumption {$({\ddot\Psi})$} we have $d_n:=\log\max_{i\le
n}\psiddot
(\lambda_{i,N}) =
o_\ff(\log n)$. By triangular inequality and decomposition (\ref
{tetaeta}), we have
%
\begin{eqnarray}
\label{AAD} \|\tS\tA_n\tS- A_n\|_2 &\leq& \|
\Pi\|_2 + \biggl\|(n-1)^{-1}\sum_{i\leq
n}
\psiddot(\lambda_{i,N}) \bigl(\tS\tilde{\eta}_i\tilde{\eta
}_i'\tS- \eta_i\eta_i'
\bigr)\biggr\|_2
\nonumber
\\
&\leq& o_\ff(1) + O_\ff\bigl(n^{-1}e^{d_n}
\bigr)\biggl\|\sum_{i\leq
n}\tilde{\delta}_i
\tilde{\delta}_i'\biggr\|_2
\\
&&{} + O_\ff\bigl(n^{-1} \bigr) \biggl\|\sum
_{i\leq n} \psiddot(\lambda_{i,N}) \bigl(\tilde{
\delta}_i\eta_i'+ \eta_i
\tilde{\delta}_i'\bigr)\biggr\|_2 \qquad\mbox{on }
\txx_{\varepsilon,
n}.
\nonumber
\end{eqnarray}
Uniformly over all unit vectors $u$ in $\RR^{N+1}$, we have
\[
u' \biggl(\sum_{i\leq n}\tilde{
\delta}_i\tilde{\delta}_i' \biggr)u \leq
\sum_{i\leq n}|\tilde{\delta}_i|^2
\leq {n\max_{i\leq n}}|\tilde{\delta}_i|^2 =
O_\ff\bigl(N^{3+\alpha}\log n \bigr) \qquad\mbox{on }
\txx_{\varepsilon,n}
\]
and by the Cauchy--Schwarz inequality,
\begin{eqnarray*}
\biggl|u' \biggl(\sum_{i\leq n}\psiddot(
\lambda_{i,N}) \bigl(\tilde{\delta}_i\eta_i'
+ \eta_i\tilde{\delta}_i'\bigr)
\biggr)u\biggr| &\leq& {O_\ff\bigl(n^{1/2}e^{d_n} \bigr)
\max_{i\leq n}}|\tilde{\delta}_i|\biggl\|\sum
_{i\leq n}\eta_i\eta_i'
\biggr\|_2^{1/2}
\\
&=& O_\ff\bigl(e^{d_n}\sqrt{n\log n}N^{(3+\alpha)/2}
\bigr) \qquad\mbox{on } \txx_{\varepsilon, n}.
\end{eqnarray*}
Therefore, the following two bounds hold:
\begin{eqnarray*}
\biggl\|\sum_{i\leq n}\tilde{\delta}_i\tilde{
\delta}_i'\biggr\|_2 &=& O_\ff
\bigl(N^{3+\alpha}\log n \bigr) \qquad\mbox{on } \txx_{\varepsilon,n},
\\
\biggl\|\sum_{i\leq n}\psiddot(\lambda_{i,N}) \bigl(
\eta_i\tilde{\delta}_i' + \tilde{\delta}
\eta_i'\bigr)\biggr\|_2 &=& O_\ff
\bigl(e^{d_n}\sqrt{n\log n} N^{(3+\alpha
)/2} \bigr) \qquad\mbox{on }
\txx_{\varepsilon,n}.
\end{eqnarray*}
By plugging into (\ref{AAD}), we can obtain that $\|\tS\tA_n\tS
-A_n\|_2 = o_\ff(1)$ on $\txx_{\varepsilon,n}$.

\begin{supplement}
\stitle{Supplement to ``Estimation in functional regression for general
exponential families.''}
\slink[doi]{10.1214/12-AOS1027SUPP} 
\sdatatype{.pdf}
\sfilename{aos1027\_supp.pdf}
\sdescription{We introduce some useful results in spectral theory and
perturbation theory in general Hilbert spaces. They serve
as powerful tools that allow us to tackle some of the statistical
approximation problems in an elegant way.
Some of the results are well-established, while others we believe are new.}
\end{supplement}


\printaddresses

\end{document}